\documentclass[11pt]{article}

\usepackage{latexsym}
\usepackage{amsfonts}
\usepackage{amsmath}

\usepackage[dvips]{color}

\newtheorem{thrm}{Theorem}[section]
\newtheorem{lemma}[thrm]{Lemma}
\newtheorem{prop}[thrm]{Proposition}

\newtheorem{cor}[thrm]{Corollary}
\newtheorem{remark}[thrm]{Remark}

\numberwithin{equation}{section}

\usepackage{enumerate}

\usepackage{amssymb}
\usepackage{mathtools}
\usepackage{mathrsfs}
\usepackage{comment}
\usepackage{hyperref}
\usepackage{xcolor}

\def\E{\mathbb{E} }
\def\P{\mathbb{P} }
\def\Q{\mathbb{Q} }
\def\R{\mathbb{R} }
\def\N{\mathbb{N} }

\def\d{\mathrm{d}}

\def\bP{{\bf P} }

\makeatletter
\begin{document}
	\allowdisplaybreaks

	\title{\Large \bf{Minimum and extremal process for a branching random walk outside the boundary case}
	\author{ \bf  
		Xinxin Chen\footnote{The research of this author is supported
			by National Natural Science Foundation of China (Grant No. 12571148) and National Key R\&D Program of China (No. 2022YFA1006500)
		} \hspace{1mm}\hspace{1mm} and \hspace{1mm}\hspace{1mm}
		Haojie Hou\thanks{The research of this author is supported by the China	Postdoctoral Science Foundation (No. 2024M76411)}}
	\hspace{1mm} }
\date{}
\maketitle

\begin{abstract}
	This work extends the studies on the minimum and extremal process of a supercritical branching random walk outside the boundary case which cannot be reduced to the boundary case. We study here the situation where the log-generating function explodes at $1$ and the random walk associated to the spine possesses a stretched exponential tail with exponent $b\in(0,\frac12)$. Under suitable conditions, we confirm the conjecture of Barral, Hu and Madaule [Bernoulli 24(2) 2018 801-841], and obtain the weak convergence for the minimum and the extremal process. We also establish an a.s. infimum result over all infinity rays of this system.
\end{abstract}

\medskip


\medskip


\section{Introduction and main results}

\subsection{Background introduction}

A branching random walk on the real line $\mathbb{R}$ is a discrete-time Markov process defined as follows. Initially, at generation $n=0$, a single particle is located at the origin. At generation $n=1$, this particle dies and produces a random number of offspring, whose positions are given by an i.i.d. copy of a point process $\mathscr{L}$. In generation $n=2$, each particle alive at generation $1$ independently reproduces in the same manner as its parent. Specifically, a particle at position $x$ generates offspring whose positions are determined by an independent copy of the translated point process $x + \mathscr{L}$. The process then continues inductively in this fashion.

\noindent We denote by \(\mathbb{P}\) the probability measure associated with the branching random walk. Let \(\mathbb{T}\) be the Galton-Watson tree describing the genealogy of the process, with root denoted by \(\varnothing\). For a particle \(u \in \mathbb{T}\), let \(|u| \in \mathbb{N}_0 := \{0,1,\dots\}\) denote its generation and \(V(u) \in \mathbb{R}\) its spatial position. We assume that the branching random walk is supercritical, i.e. $\E(\#\mathscr{L})>1$, so that it survives with positive probability. We also assume that $\E(\#\mathscr{L})<\infty$ \footnote{If $\E(\#\mathscr{L})=\infty$, the system may grow double exponentially.}so that the tree $\mathbb{T}$ grows exponentially. For each \(\beta \in \mathbb{R}\), we define the log-generating function
\[
\phi(\beta):= \log \, \mathbb{E}\!\left[ \sum_{|u|=1} e^{-\beta V(u)} \right]=\log \E\left[\int e^{-\beta x}\mathscr{L}(\d x)\right] \; \in \; (-\infty, \infty].
\]
We also assume that the point process \(\mathscr{L}\) is not almost surely supported on a deterministic lattice. Then $\phi$ is strictly convex on its domain $\text{dom}(\phi)=\{\beta: \phi(\beta)<\infty\}$ as long as $\text{dom}(\phi)$ is non-trivial.

In this work, we are interested in the minimal position of the system at generation $n$, which is defined by 
\[
M_n:=\inf_{|u|=n}V(u)
\]
with the convention that $\inf\emptyset =  \infty$.

\noindent It is well known (see, e.g., \cite{Biggins1976, Hammersley1974, Kingman1975}) that if $\{\beta>0: \phi(\beta)<\infty\}\neq \emptyset$, then conditioned on survival $\{\mathbb{T} = \infty\}$, a.s.,
\[
\lim_{n\to\infty} \frac{M_n}{n}=v:= -\inf_{\beta>0}\frac{\phi(\beta)}{\beta}\in\mathbb{R}.
\]
Hammersley \cite{Hammersley1974} raised the question of the second order of $M_n$, that is, the asymptotic of $M_n-vn$. This is deeply related to the way that $\inf_{\beta>0}\frac{\phi(\beta)}{\beta}$ is achieved. There exist different cases:
\begin{itemize}
	\item[(I)] $\inf_{\beta>0}\frac{\phi(\beta)}{\beta}$ is attained at $\beta_v\in(0,\infty)$, $(0,\beta_v]\subset \text{dom}(\phi)$ and $\frac{\phi(\beta_v)}{\beta_v} = \phi'(\beta_v-)$;
	\item[(II)] $\inf_{\beta>0}\frac{\phi(\beta)}{\beta}$ is attained at $\beta_v\in(0,\infty)$, $(0,\beta_v]\subset \text{dom}(\phi)$ and $\frac{\phi(\beta_v)}{\beta_v} >\phi'(\beta_v-)$ and $\phi(\beta)=\infty, \forall \beta>\beta_v$;
	\item[(III)] $\inf_{\beta>0}\frac{\phi(\beta)}{\beta}=\lim_{\beta\to\infty}\frac{\phi(\beta)}{\beta}$ and $(0,\infty)\subset \text{dom}(\phi)$.
\end{itemize}
In the first two cases (I)-(II), we can make a linear transformation $(u, V(u))\mapsto (u, \widetilde{V}(u):=\beta_vV(u)+\phi(\beta_v)|u|)$ and see that $\frac{\inf_{|u|=n}\widetilde{V}(u)}{n}\to v=0$ a.s. on $\{\mathbb{T} = \infty\}$. We thus could set that $\phi(1)=\phi'(1-)=0$ for the case (I), and set that $\phi(1)=0$,$\phi'(1-)<0$ and $\phi(\beta)=\infty,\forall\beta>1$ for the case (II).

\noindent Following \cite{BK2005} and \cite{BHM2018}, we say that the case (I) with $\phi(1)=\phi'(1-)=0$ is the boundary case, whereas the case (II) and the case (II) are both ``outside the boundary case".

\paragraph{Boundary case (I).} For the boundary case, the resolution of Hammersley's problem was advanced by several significant contributions (see, e.g., \cite{AR2009, HuShi2009} and references therein), culminating in the sharp result obtained by Aïdékon \cite{Aidekon2013}, which shows that under mild moment conditions, there exist positive constants $\lambda, C_0$ such that 
\[
\lim_{n\to\infty}\P(M_n > \lambda \log n +x) = \E\left[ \exp\{-C_0 e^x D_\infty\}\right], \forall x\in \R,
\]
where $D_\infty\ge0$ is a.s. limit of the so-called derivative martingale. Further, the extremal point process $\sum_{|u|=n}\delta_{V(u) -\lambda \log n}$ is investigated by \cite{Madaule}\footnote{see \cite{ABBS, ABK} for the analog for branching Brownian motion} and it turns out that in the vague topology, $\sum_{|u|=n}\delta_{V(u) -\lambda \log n}$ converges in law to some Decorated Poisson point process.

Let us discuss also some asymptotical behaviors along the infinite rays. We say that $\Theta:= \{\Theta_0=\emptyset, \Theta_1,...,\Theta_,...\}\subset \mathbb{T}$ is an infinite ray of $\mathbb{T}$ if for all $n\geq 0$, $\Theta_n$ is the parent of $\Theta_{n+1}$. Define $\partial \mathbb{T}$ to be the set of all the infinite rays. In the boundary case with some further integrability conditions, according to Jaffuel \cite{Jaffuel}, there exists an explicit constant $a_0>0$ such that on  $\{\mathbb{T} = \infty\}$, a.s., 
\[
\inf_{\Theta \in\partial\mathbb{T}}\limsup_{n\to\infty} \frac{V(\Theta_n)}{n^{1/3}} = a_0.
\]

\paragraph{Outside boundary case: case (III).} One can refer to Bramson \cite{Bramson} for an example in the case (III), where $M_n-vn$ is of order $\log\log n$. It is assumed in  \cite{Bramson} that the displacements are bounded and 
\[
\E\bigg[\sum_{|u|=1}1_{V(u) = \textrm{ess}\inf\mathscr{L}}\bigg]=1.
\]

\paragraph{Outside boundary case: case (II)} In this work, we study the minimum of branching random walks within the non-boundary regime, focusing on a specific scenario stated below. First, throughout the paper, we assume that
\begin{equation}\label{Assumption1}
	\phi (1)=0,\quad \phi'(1-)\in (-\infty, 0) \quad\mbox{and}\quad \phi(1+)=\infty. 
\end{equation}
Next, it is convenient to frame our setting in terms of the law of $X$ under the probability measure $\bP$, which is defined by 
\begin{equation}\label{Law-of-spine}
	\mathbf{E}[f(X)]=\int_\R f(x) \bP(X\in\d x) := \E\left[\sum_{|u|=1}f(V(u))e^{-V(u)} \right],
\end{equation}
for any bounded measurable function $f$. Note immediately that $\phi(\beta) = \log \mathbf{E}[e^{-(\beta-1)X}]$. The property $\phi(1+)=\infty$ then follows, in particular, when $X$ has a density with polynomial or sub-exponential decay in the left tail. 

\noindent The present work addresses the regime of subexponential decay with exponent $b\in(0,\frac12)$, in contrast to the polynomial decay case treated by Barral, Hu, and Madaule \cite{BHM2018}. More precisely, in \cite{BHM2018}, under \eqref{Assumption1}, the so-called LlogL condition (see \eqref{Assumption2} below) and the condition that there exist some constants $\gamma>3,\alpha>2, x_0<0$ and a slowly varying function $\ell$ on $-\infty$ such that 
\begin{align}\label{Assumption3a}
	\mathbf{E}\left( (\max\{X, 0\})^\gamma \right)<\infty\quad \mbox{and}\quad \mathbf{P}(X\leq x)=\int_{-\infty}^x \ell(y)|y|^{-\alpha} \mathrm{d}y,\quad\forall x\leq x_0,
\end{align}
it is proved that 
\begin{align}\label{Thm-Poly-tail}
	\lim_{n\to\infty} \P(M_n\geq \alpha\log n -\log \ell(-n) +x )= \E\left(\exp\left\{-m^{-\alpha}C^* e^x W_\infty\right\}\right), \forall x\in \R,
\end{align}
where $W_\infty$ is the a.s. limit of the additive martingale $W_n=\sum_{|u|=n}e^{-V(u)}$ and
\begin{align}\label{Def-of-C-star}
	C^*:= \sum_{j=0}^\infty \E(e^{-M_j})\in(0,\infty). 
\end{align}
Note that, in \cite[Remark 1.6]{BHM2018}, the authors conjectured that  if \eqref{Assumption3a} is replaced by 
\begin{align}\label{Assumption3b}
	\mathbf{P}(X\leq x)=\int_{-\infty}^x \ell(y)|y|^a e^{-\lambda |y|^b}\mathrm{d}y,\quad\forall x\leq x_0,
\end{align}
with $a\in \R, \lambda>0, b\in (0,1)$ and $\ell$ a non-negative function such that $\lim_{x\to-\infty} \ell(x)=: \ell_\infty\in (0,\infty)$, then $M_n$ should be of order $n^b$. Furthermore, if $b\in (0, \tfrac{1}{2})$, then $M_n- \lambda(mn)^b+a\log n$ is tight. In this paper, we confirm this conjecture for $b\in (0, \tfrac{1}{2})$, by establishing the weak convergence of the minimum $M_n$ and the associated extremal process.

Concerning the asymptotical behaviors along the infinite rays, recently, it is proved by A\"id\'ekon,  Hu and Shi \cite[Theorem 1.2]{AHS} that  under a slightly general assumption than \eqref{Assumption3a}(corresponding to \cite{BHM2018}), there exists a positive constant $a_1$ such that conditioned on survival, a.s.,
\begin{align}\label{Limit-of-Ray}
	\inf_{\Theta\in \partial\mathbb{T}}\limsup_{n\to\infty} \frac{V(\Theta_n)}{\sqrt{n\log n}} = a_1.
\end{align}
It is evident that the asymptotic phenomenology differs significantly between the boundary and non-boundary regimes.

\subsection{Main results}
Now we are ready to state our setting. To simplify the life, we suppose that the point process $\mathscr{L}$ is of the following form:
\begin{align}\label{IID-case}
	\mathscr{L}:= \sum_{i=1}^\nu \delta_{Y_i},
\end{align}
where $\nu\in\N_0$ is the offspring law of GW tree $\mathbb{T}$ and  $\{Y_i\}_{i\in \mathbb{N}}$ is a family of i.i.d. real-valued random variables independent of $\nu$. We assume that  \eqref{Assumption1} holds for $\mathscr{L}$.

\noindent Recall the law of $X$ given in \eqref{Law-of-spine}. We assume that \eqref{Assumption3b} holds for $X$. And \eqref{Assumption1} shows that
\begin{align}\label{Def-of-alpha}
	m:= -\phi'(1-)= \mathbf{E}(X)\in (0,\infty).
\end{align}

\noindent \eqref{Assumption1} or in fact $\phi(1)=0$, implies also that $W_n:=\sum_{|u|=n}e^{-V(u)}$ is a martingale with respect to the natural filtration $\mathcal{F}_n:=\sigma((u,V(u)); |u|\le n)$. It is immediate that $W_n$ converges a.s. to some limit $W_\infty\ge0$. In addition,  it is proved in \cite{Biggins1976, DL1983, KaPe1976} that  $W_n$ converges in $L^1(\P)$ to $W_\infty$ if and only if the so-called LlogL condition 
\begin{align}\label{Assumption2}
	\E\left( W_1 \log_+ W_1\right)<\infty,
\end{align}
holds\footnote{where $\log_+x:= \max\{ \log x, 0 \}$}. Note that, under \eqref{Assumption2}, a.s., $\{W_\infty>0\}=\{\# \mathbb{T}=\infty\}$. We assume that \eqref{Assumption2} holds in the following.

Recall $C^*$ in \eqref{Def-of-C-star}. Further, we can define
\begin{align}\label{Def-of-C-star-f}
	C^*(f): = \ell_\infty m^a    \sum_{j=0}^\infty \E  \left[  e^{-M_j}   \left(1+ \int_0^\infty e^{z} \left(1- e^{-\sum_{|v|=j} f(V(v)-M_j+z)} \right) \mathrm{d} z\right)\right] ,
\end{align}
for any $f\in \mathcal{S}$ where 
\begin{align}\label{Def-of-test-function}
	\mathcal{S}:= \left\{f: \ f\ \mbox{is a non-negative continuous function with}\ \mbox{supp}(f)\subset (-\infty, R_f)\ \mbox{for some } R_f>0\right\}.
\end{align}
In particular,  $C^*(0)= \ell_\infty m^a C^*$. We will show that $C^*(f)\in(0, \infty)$ in our setting.

Set 
\begin{align}\label{Def-of-alpha}
	\alpha_n:= \lambda(mn)^b- a\log n \quad\mbox{with}\quad m= \mathbf{E}(X)\in (0,\infty).
\end{align}

Our first main result is stated as follows. 

\begin{thrm}\label{thm1}
	Assume \eqref{Assumption1}, \eqref{Assumption3b} and \eqref{Assumption2} with $b\in (0, \frac{1}{2})$. Then for any $x\in \R$ and $f\in \mathcal{S}$, 
	\begin{align}\label{Joint-prob}
		\lim_{n\to\infty}  \E \left[ e^{-\sum_{|u|=n} f(V(u)- \alpha_n )} 1_{\left\{ M_n> \alpha_n + x \right\}}\right]=  \E\left[ \exp\left\{- C^*(f_x)e^{x}W_\infty \right\}\right],
	\end{align}
	where $C^*(f_x)$ is given as in \eqref{Def-of-C-star-f} with $f_x:=f(\cdot+x)$. Moreover, for any non-negative function $f$ with bounded support, we have
	\begin{align}\label{Extremal-process}
		\lim_{n\to\infty}  \E \left[ e^{-\sum_{|u|=n} f(V(u)- \alpha_n )}\right] = \E\left[ \exp\left\{- \left(C^*(f)-C^*(0)\right)W_\infty \right\}\right]. 
	\end{align} 
\end{thrm}

\begin{remark}
	In the appendix of \cite{BHM2018}, the authors briefly addressed the weak convergence of $M_n$ under more general conditions on $\mathscr{L}$. We note that within our framework, assumptions analogous to \cite[(A.1)--(A.3)]{BHM2018} remain applicable. Consequently, following a parallel argument, Theorem \ref{thm1} holds with a modified functional $C^*(f)$; the detailed verification is omitted here.
\end{remark}
Taking $f=0$ in Theorem \ref{thm1}, we immediately get the following result.
\begin{cor}
	Under the same assumptions as Theorem \ref{thm1}, for any $x\in \R$, 
	\begin{align}
		\lim_{n\to\infty}  \P \left(M_n> \alpha_n + x \right) =  \E\left( \exp\left\{- C^*e^{x}W_\infty \right\}\right).
	\end{align}
\end{cor}

Next, we consider the weak convergence of the extremal point process defined by
\[
\mathcal{E}_n:= \sum_{|u|=n} \delta_{V(u)-\alpha_n}.
\]

Our second result is stated as follows.
\begin{thrm}\label{thm2}
	Assume \eqref{Assumption1}, \eqref{Assumption3b} and \eqref{Assumption2} with $b\in (0, \frac{1}{2})$.  $\mathcal{E}_n$ converges in distribution to $\mathcal{E}_\infty$ in the sense of vague topology. The limiting extremal process $\mathcal{E}_\infty$ is defined by
	\[
	\mathcal{E}_\infty:= \sum_{i=1}^\infty 1_{\{M_{q_i}^{(i)}\geq p_i\}} \sum_{u\in \mathbb{T}^{(i)}: |u|=q_i} \delta_{V(u)- p_i} ,
	\]
	where 
	\begin{itemize}
		\item given $W_\infty$,  $\mathcal{P}:= \sum_{i=1}^\infty \delta_{(p_i, q_i)}$ is a Poisson point process with intensity $W_\infty \ell_\infty m^a e^{-z}\mathrm{d}z \otimes \delta_{\mathbb{N}_0}(\mathrm{d} n)$;
		\item for every $i\ge 1$, $\{M_n^{(i)}: n \in \mathbb{N}_0\}$ is the minimal position process of the branching random walk $\{V^{(i)}(u), u\in \mathbb{T}^{(i)}\}$ with $\{V^{(i)}(u), u\in \mathbb{T}^{(i)}\}$ being i.i.d. copies of the branching random walk $\{V(u), u\in \mathbb{T}\}$ and being independent of $\mathcal{P}$.
	\end{itemize}
	
\end{thrm}

Our final result concerns the asymptotical behaviors of the infinite rays, complementing the result of \cite[Theorem 1.2]{AHS} in our case. 
\begin{thrm}\label{thm3}
	Assume \eqref{Assumption1}, \eqref{Assumption3b} and \eqref{Assumption2} with $b\in (0, \frac{1}{2})$.  There exists a positive constant $a_*$ such that on $\{\# \mathbb{T}=\infty\}$, a.s.,
	\[
	\inf_{\Theta\in \partial\mathbb{T}}\limsup_{n\to\infty} \frac{V(\Theta_n)}{n^{\frac{1}{2-b}}} = a_*.
	\]
\end{thrm}

\subsection{Proof strategies and discussions}

Our proof strategy for Theorems \ref{thm1} and \ref{thm2} adapts the framework of \cite{BHM2018}, originally developed by \cite{Aidekon2013}. The key mechanism driving the weak convergence of $M_n$ is that, with high probability, there is exactly one large jump whose timing is near $n$. In contrast to the polynomial decay case \cite{BHM2018}, where the jump size scales as $-\frac{n}{(\log n)^3} = -n^{1+o(1)}$, the stretched exponential setting necessitates a finer tuning. We ultimately set the jump size to be $-(mn - A n^{1-b}\log n)$ for appropriate constants $m$ and $A$.

\noindent We briefly explain here why the regime \( b \in [\tfrac12, 1) \) is not treated in this work. 
Let \( \xi_n \) be a particle attaining the minimum at generation \( n \), i.e., \( V(\xi_n) = M_n \), 
and let \( T \le n \) be the time at which a large jump occurs, with jump size  
\( \zeta := V(\xi_T) - V(\xi_{T-1}) \).

\noindent On the one hand, the second‑order fluctuation of the random walk \( V(\xi_n) - \zeta \) is of order \( n^{1/2} \). On the other hand, for \( b \in (\tfrac12, 1) \) we have \( n^{1-b}\log n = o(n^{1/2}) \). 
Consequently, the natural scaling for the large jump would become \( -(mn - O(\sqrt{n})) \) 
instead of \( -(mn - A n^{1-b}\log n) \). This shift shows that the density around the 
large jump \( \zeta \) becomes strongly coupled with the Gaussian fluctuations of the 
random walk, thereby changing the nature of the problem.

\noindent The borderline case \( b = \tfrac12 \) requires a different and more delicate analysis 
compared to \( b \in (0, \tfrac12) \), and is therefore left for future investigation.

Our proof strategy for Theorem~\ref{thm3} closely mirrors that developed in \cite{AHS}. It rests on three main steps: (i) the branching property implies that the limit \(a_*\) is deterministic; (ii) the positivity \(a_*>0\) is derived from a first-moment method; and (iii) the finiteness \(a_*<\infty\) is obtained via a coupling argument.

The rest of the paper is organized as follows. Section~\ref{S2} establishes elementary properties of the random walk $S_n$. The proofs of Theorems~\ref{thm1} and~\ref{thm2} are presented in Section~\ref{S3}. Section~\ref{S5} contains the proof of Theorem~\ref{thm3}. Finally, the proof of Proposition~\ref{prop:key-proposition} is given in Section~\ref{S4}.

Notation convention: We use $a_n \lesssim b_n$ or $a_n= O(b_n)$ to mean that that there exist constants $N, C$ such that for any $n\geq N$, $a_n\leq C b_n$. $a_n=o(1)$ means that $\lim_{n\to\infty} a_n =0$. We also use $a_n \lesssim_{K} b_n (a_n \lesssim_{K, L} b_n)$ to denote that  there exist constants $N= N(K), C=C(K) (N= N(K,L), C=C(K,L))$ such that for any $n\geq N$, $a_n\leq C b_n$. Denote by $a_n\gtrsim b_n (a_n\gtrsim_K b_n, a_n\gtrsim_{K,L} b_n )$ if $b_n\lesssim a_n (b_n\lesssim_K a_n, b_n\lesssim_{K,L} a_n )$. Denote by $a_n \asymp b_n \left(a_n \asymp_K b_n, a_n \asymp_{K,L} b_n\right)$ if $a_n\gtrsim b_n (a_n\gtrsim_K b_n, a_n\gtrsim_{K,L} b_n )$ and $a_n\lesssim b_n (a_n\lesssim_K b_n, a_n\lesssim_{K,L} b_n )$.
Since in the whole paper, we regard $\lambda, m, b, a, \ell_\infty$ as known constants, so all the constants may depend on these parameters.

\section{Elementary properties for random walk $S_n$}\label{S2}

Let $\{S_n, n\in \mathbb{N}_0,  \mathbf{P}\}$ be a random walk with $S_0= 0$ such that $(S_{n}- S_{n-1})_{n\geq 1}$ are i.i.d. copies of $(X, \mathbf{P})$ given as in \eqref{Law-of-spine}. We state some preliminary results on this random walk.

Under the assumption \eqref{IID-case}, our conditions \eqref{Assumption1}, \eqref{Assumption3b} and \eqref{Assumption2} are equivalent to
\begin{equation}\label{Equal-Conditions}
	\begin{split}
		&\E(\nu)\in (1,\infty), \quad \E(\nu \log_+ \nu)<\infty,\quad  \E(e^{-Y_1})= \frac{1}{\E (\nu)},\quad  \E(Y_1e^{-Y_1})>0 ,\\
		& \P(Y_1\leq x) =\frac{1}{\E (\nu)}\int_{-\infty}^x \ell(y)|y|^a e^{-\lambda |y|^b +y}\mathrm{d}y,\quad\forall x\leq x_0 . 
	\end{split}
\end{equation}
Therefore, it is easy to see from \eqref{Equal-Conditions} that 
\begin{align}\label{Assumption4}
	\mathbf{E}(|X|^k) =\E(\nu) \E(|Y_1|^ke^{-Y_1}) <\infty,\quad \mbox{for all } k\in \N.
\end{align}
For each $n\in \N$, define
\begin{align}\label{Def-of-zeta}
	\zeta_n:= mn - A_1 n^{1-b} \log n,\quad \widehat{\zeta}_n:= \zeta_n+m\quad\mbox{and}\quad \theta_n:= \lambda \widehat{\zeta}_n^{b-1}- \frac{A_2 \log \widehat{\zeta}_n}{\widehat{\zeta}_n},
\end{align}
where $A_1$ and $A_2$ are two fixed constants such that 
\begin{align}\label{Choice-of-A}
	(1-b)\lambda  m^{b-1} A_1> A_2-a + 2 \quad \mbox{and}\quad A_2> a+1+2(1-b).
\end{align}
We always assume that $n$ is large enough such that $\theta_n>0$. Define the centralized r.v.
\begin{align}
	\widehat{X}:=  X-m.
\end{align}
Since $\lim_{y\to-\infty} |y|/|y-m| =1$ and $ \lim_{y\to-\infty} (|y|^b-|y-m|^b )=0$, by \eqref{Assumption3b}, 
the function $\widehat{\ell}$ defined by
\begin{align}\label{Density-hat-X}
	\int_{-\infty}^x \widehat{\ell}(y)|y|^a e^{-\lambda |y|^b}\mathrm{d}y:=\mathbf{P}(\widehat{X}\leq x)=\int_{-\infty}^x \ell(y-m)|y-m|^a e^{-\lambda |y-m|^b}\mathrm{d}y, \quad\forall x\leq x_0-m
\end{align}
satisfies $\lim_{y\to-\infty} \widehat{\ell}(y)=\ell_\infty$. According to elementary calculation, it holds that
\begin{align}\label{Step1}
	\mathbf{P}(\widehat{X}\leq -\widehat{\zeta}_n)&=\mathbf{P}(X\leq -\zeta_n) \asymp \int_{\zeta_n}^\infty  y^a e^{-\lambda y^b}\mathrm{d}y \stackrel{z:= \lambda y^b}{=} \frac{1}{b \lambda^{(a+1)/b}} \int_{\lambda \zeta_n^b}^\infty  z^{(a+1-b)/b} e^{-z}\mathrm{d}z\nonumber\\
	& \stackrel{x:=z-\lambda \zeta_n^b}{=} \frac{1}{b\lambda} \zeta_n^{a+1-b} e^{-\lambda \zeta_n^b}\int_{0}^\infty  \left(1+ \frac{x}{\lambda \zeta_n^b}\right)^{(a+1-b)/b} e^{-x}\mathrm{d}x\nonumber\\
	&\asymp \zeta_n^{a+1-b} e^{-\lambda \zeta_n^b},
\end{align}
where in the last inequality we used the fact that $\lim_{n\to\infty} \int_{0}^\infty  \left(1+ \frac{x}{\lambda \zeta_n^b}\right)^{(a+1-b)/b} e^{-x}\mathrm{d}x=1$ according to dominated convergence theorem. 
Moreover, combining Taylor's expansion and \eqref{Step1},
\begin{align}\label{Step16}
	\mathbf{P}(X\leq -\zeta_n) & \asymp n^{a+1-b}e^{-\lambda (mn)^b + \lambda b (mn)^{b-1} A_1 n^{1-b} \log n +o(1)}\nonumber\\
	&\asymp e^{-\alpha_n} n^{1-b+ \lambda b m^{b-1} A_1}. 
\end{align}

Define 
\begin{align}\label{Def-of-large-jump}
	\tau_{\zeta}:= \min\left\{k: X_k <-\zeta \right\}\quad\mbox{and}\quad \tau_{\zeta}^{(2)}:= \min\left\{k> \tau_{\zeta}: X_k <-\zeta\right\}.
\end{align}

\begin{lemma}\label{lem1}
	\begin{itemize}
		\item[(i)] 
		Assume \eqref{Assumption3b} holds with  $b\in (0,1)$. Then for large $n$, we have
		\[
		\mathbf{E}\left(\exp\left\{- \theta_n \max\left\{\widehat{X}, -\widehat{\zeta}_n \right\}  \right\}\right) -1 \lesssim \theta_n^2. 
		\]
		
		\item[(ii)] Assume \eqref{Assumption3b} holds with $b\in (0,\tfrac{1}{2})$. Then, when $n$ is large enough, for all $x\geq 0$ and $1\leq q\leq 2n$, 
		\begin{align}
			\mathbf{P}\left(S_q - m q \leq -x, \tau_{\zeta_n} >q\right) \lesssim e^{-\theta_n x}.
		\end{align}
	\end{itemize}
\end{lemma}
\textbf{Proof: } 
(i)
From \eqref{Choice-of-A}, there exists a sufficiently small constant $\varepsilon>0$ such that $(1-\varepsilon)A_2\geq  a+1+ 2(1-b)$. Let $n$ be large enough such that $\theta_n \widehat{\zeta}_n= \lambda\widehat{\zeta}_n^b- A_2\log \widehat{\zeta}_n >\frac{1}{1-\varepsilon}$. First noticing that 
\begin{align}\label{Step2}
	&\mathbf{E}\left(\exp\left\{- \theta_n \max\left\{\widehat{X}, -\widehat{\zeta}_n \right\}  \right\}\right) -1 \nonumber\\
	&= e^{\theta_n \widehat{\zeta}_n} \mathbf{P}(\widehat{X}\leq -\widehat{\zeta}_n)+ \left(\mathbf{E}\left(e^{-\theta_n \widehat{X}} 1_{\{\theta_n \widehat{X} \geq -1 \}}\right) -1 \right)+  \mathbf{E}\left(e^{-\theta_n \widehat{X}} 1_{\{ -\widehat{\zeta}_n< \widehat{X} < -\theta_n^{-1} \}}\right)\nonumber\\
	&=: I_1+I_2+I_3.
\end{align}
For $I_1$, combining \eqref{Def-of-zeta} and   \eqref{Step1}, we have
\begin{align}\label{Step3}
	I_1\asymp  e^{\theta_n \widehat{\zeta}_n} \zeta_n^{a+1-b} e^{-\lambda \zeta_n^b}  \asymp e^{\theta_n \widehat{\zeta}_n} \widehat{\zeta}_n^{a+1-b} e^{-\lambda \widehat{\zeta}_n^b} = \widehat{\zeta}_n^{a+1-b} e^{-A_2\log \widehat{\zeta}_n}< \widehat{\zeta}_n^{-b-2(1-b)}\lesssim \theta_n^2,
\end{align}
where the last inequality follows from $\theta_n \asymp \widehat{\zeta}_n^{b-1}$.  For $I_2$, noticing that $|e^{-x}-1+x|\lesssim x^2$ for any $x\geq -1$ and that 
$\mathbf{E}\left(\theta_n \widehat{X} 1_{\{\theta_n \widehat{X} \geq -1 \}}\right) = - \mathbf{E}\left(\theta_n \widehat{X} 1_{\{\theta_n \widehat{X} < -1 \}}\right)>0$ since $\mathbf{E}(\widehat{X})=0$, we see that 
\begin{align}\label{Step4}
	I_2 & = \mathbf{E}\left(\left(e^{-\theta_n \widehat{X}}-1+\theta_n \widehat{X}\right) 1_{\{\theta_n \widehat{X} \geq -1 \}}\right) - \mathbf{P}(\theta_n \widehat{X}<-1)- \mathbf{E}\left(\theta_n \widehat{X} 1_{\{\theta_n \widehat{X} \geq -1 \}}\right)\nonumber\\
	& \leq  \mathbf{E}\left(\left(e^{-\theta_n \widehat{X}}-1+\theta_n \widehat{X}\right) 1_{\{\theta_n \widehat{X} \geq -1 \}}\right)\nonumber\\
	&\lesssim \theta_n^2 \mathbf{E}(\widehat{X}^2)\lesssim \theta_n^2.
\end{align}
Now we treat $I_3$. It follows from \eqref{Density-hat-X} that
\begin{align}\label{Step5}
	I_3 & \lesssim \int_{\theta_n^{-1}}^{\widehat{\zeta}_n} y^a e^{-\lambda y^b} e^{\theta_n y}\mathrm{d} y= \left( \int_{\theta_n^{-1}}^{(1-\varepsilon)\widehat{\zeta}_n} +  \int_{(1-\varepsilon)\widehat{\zeta}_n}^{\widehat{\zeta}_n}\right)y^a e^{-\lambda y^b} e^{\theta_n y}\mathrm{d} y.
\end{align}
According to the definition of $\theta_n$ in \eqref{Def-of-zeta}, for all $0<y\leq \widehat{\zeta}_n$,
\begin{align}\label{Step6}
	\theta_n y-\lambda y^b &= - \frac{A_2 \log \widehat{\zeta}_n}{\widehat{\zeta}_n} y + \lambda y^b \left(\widehat{\zeta}_n^{b-1} y^{1-b}-1\right)\nonumber\\
	& \leq \begin{cases}
		- \frac{A_2 \log \widehat{\zeta}_n}{\widehat{\zeta}_n} y ,\quad & (1-\varepsilon)\widehat{\zeta}_n\leq y\leq \widehat{\zeta}_n;\\
		-\lambda (1-(1-\varepsilon)^{1-b})y^b,\quad &y\leq (1-\varepsilon)\widehat{\zeta}_n.
	\end{cases}
\end{align}
Plugging \eqref{Step6} into \eqref{Step5} yields that
\begin{align}\label{Step7}
	I_3& \lesssim \int_{\theta_n^{-1}}^{(1-\varepsilon)\widehat{\zeta}_n} y^a e^{-\lambda (1-(1-\varepsilon)^{1-b}) y^b}\mathrm{d} y+ \int_{(1-\varepsilon)\widehat{\zeta}_n}^{\widehat{\zeta}_n}y^a  e^{-\frac{A_2 \log \widehat{\zeta}_n}{\widehat{\zeta}_n} y}\mathrm{d} y\nonumber\\
	& \leq \theta_n^2 \int_{\theta_n^{-1}}^{\infty} y^{a+2} e^{-\lambda (1-(1-\varepsilon)^{1-b}) y^b}\mathrm{d} y + e^{-A_2 (1-\varepsilon) \log \widehat{\zeta}_n}\int_{(1-\varepsilon)\widehat{\zeta}_n}^{\widehat{\zeta}_n}y^a  \mathrm{d} y\nonumber\\
	& \lesssim \theta_n^2 + \widehat{\zeta}_n^{-2(1-b)}\lesssim \theta_n^2. 
\end{align}
Combining \eqref{Step2}, \eqref{Step3}, \eqref{Step4} and \eqref{Step7}, we complete the proof of (i). 

(ii) 
Define $\widehat{X}_i:= X_i-m$ and $\widehat{S}_q^{(n)}:= \sum_{i=1}^q \max\left\{\widehat{X}_i, -\widehat{\zeta}_n \right\}$, then by Markov's inequality, 
\begin{align}\label{Step8}
	\mathbf{P}\left(S_q - m q \leq -x, \tau_{\zeta_n} >q\right)  &  = 	\mathbf{P}\Big(S_q - m q \leq -x,  \min_{1\leq j\leq q} X_j \geq -\zeta_n\Big)  \nonumber\\
	&\leq 	\mathbf{P}\left(-\theta_n \widehat{S}_q^{(n)} \geq \theta_n x\right)\nonumber\\
	&  \leq e^{-\theta_n x} \left(\mathbf{E}\left(\exp\left\{- \theta_n \max\left\{\widehat{X}, -\widehat{\zeta}_n \right\}  \right\}\right) \right)^q.
\end{align}
According to (i), there exists a constant $C>0$ such that for large $n$, 
\[
\mathbf{E}\left(\exp\left\{- \theta_n \max\left\{\widehat{X}, -\widehat{\zeta}_n \right\}  \right\}\right)\leq 1+ C\theta_n^2 \leq e^{C\theta_n^2}.
\]
Plugging the above inequality to \eqref{Step8}, we see that for all $1\leq q\leq 2n$, 
\begin{align}
	\mathbf{P}\left(S_q - m q \leq -x, \tau_{\zeta_n} >q\right)  \leq e^{-\theta_n x} e^{C\theta_n^2 q }\leq e^{-\theta_n x} e^{C2\theta_n^2 n }.
\end{align}
Together with the fact that  $\theta_n^2 n \lesssim n^{2b-1}=o(1)$ for $b\in (0,\tfrac{1}{2})$, we complete the proof of (ii). 

\hfill$\Box$

\begin{lemma}\label{lem3}
	Assume \eqref{Assumption3b}  holds with $b\in (0,\tfrac{1}{2})$.
	
	(i) When $n$ is large enough,  we have
	\begin{align}
		\mathbf{P}\left(\tau_{\zeta_n}^{(2)}\leq n\right)\lesssim n^2 \zeta_n^{2(a+1-b)} e^{-2\lambda\zeta_n^b} . 
	\end{align}
	
	(ii) When $n$ is large enough, for any $y< 4p^{(3-2b)/4}, z\in [0, 4p^{(3-2b)/4})$, $p\in [\sqrt{n}, 2n]$ and any $1\leq i\leq p$, it holds that
	\begin{align}
		& \mathbf{P}\left( S_p-y\in [z, z+1], \min_{i\leq j\leq p} S_j\geq y, \tau_{\zeta_n}=i, \tau_{\zeta_n}^{(2)}> p\right)\nonumber\\
		&\lesssim \mathbf{P}\left(\max_{0\leq j\leq p-i} S_j \leq z+1\right) e^{-\alpha_p} +\frac{1}{n^2}e^{-\alpha_n}.
	\end{align}
\end{lemma}
\textbf{Proof: } (i) On $\{\tau_{\zeta_n}^{(2)}\leq n\}$, there exist at least two large jumps up to time $n$, which implies that 
\begin{align}\label{Step9}
	&	\mathbf{P}\left(\tau_{\zeta_n}^{(2)}\leq n\right) \leq \sum_{1\leq i\neq j \leq n} 	\mathbf{P}\left(X_i, X_j <-\zeta_n\right)\leq n^2 \mathbf{P}(X<-\zeta_n)^2.
\end{align}
Combining \eqref{Step1} and \eqref{Step9}, we get (i).

(ii) For each $1\leq i\leq p$, define 
\begin{align}\label{Def-of-E-i-p}
	E_i^{(p)}:= \left\{\left|S_p-X_i-m(p-1)\right|\leq p^{(3-2b)/4} \right\}.
\end{align}
On one hand, according to the independence between $E_i^{(p)}$ and $X_i$, 
\begin{align}\label{Step13}
	& \mathbf{P}\left( S_p-y\in [z, z+1], \min_{i\leq j\leq p} S_j\geq y, \tau_{\zeta_n}=i, \tau_{\zeta_n}^{(2)}> p, E_i^{(p)}\right)\nonumber\\
	& \leq  \mathbf{P}\left( S_p-  \min_{i\leq j\leq p} S_j \leq z+1, S_p-y\in [z, z+1], E_i^{(p)}\right)\nonumber\\
	& = \mathbf{E}\Big(1_{\left\{ S_p-  \min_{i\leq j\leq p} S_j \leq z+1 \right\}} 1_{E_i^{(p)}} \mathbf{P}\left(X_i\in[z+y-t,z+1+y-t]\right)\bigg|_{t=S_p-X_i}\Big).
\end{align}
Noticing that on $E_i^{(p)}$, when $n$ is large enough, for all $z, y< 4p^{(3-2b)/4}$ and $p\in [\sqrt{n},2n]$,
\begin{align}
	z+1+y-t\leq 8p^{(3-2b)/4}+1- (m(p-1)-p^{(3-2b)/4})< -\left(mp- 10p^{(3-2b)/4} \right),
\end{align}
which implies that 
\begin{align}\label{Step11}
	\mathbf{P}\left(X_i\in[z+y-t,z+1+y-t]\right)& \leq \sup_{r< -\left(mp- 10p^{(3-2b)/4} \right)} \int_{r-1}^r \ell(y) |y|^a e^{-\lambda |y|^b}\mathrm{d} y\nonumber\\
	& \asymp \sup_{r> mp- 10p^{(3-2b)/4} }  r^a e^{-\lambda r^b} \nonumber\\
	&\asymp \left(mp- 10p^{(3-2b)/4} \right)^a e^{-\lambda \left(mp- 10p^{(3-2b)/4} \right)^b}.
\end{align}
According to Taylor's expansion, we have
\begin{align}\label{Step12}
	\left(mp- 10p^{(3-2b)/4} \right)^b &= (mp)^b - b(mp)^{b-1} \times 10p^{(3-2b)/4} +O(p^{b-2}\times p^{(3-2b)/2})\nonumber\\
	& = (mp)^b + o(1),
\end{align}
where in the last equallity we used the fact that $b-1+ (3-2b)/4 = (2b-1)/4<0$ and that $b-2+(3-2b)/2 = -1/2<0$. Combining \eqref{Step13}, \eqref{Step11} and \eqref{Step12}, we finally conclude that 
\begin{align}\label{Step14}
	& \mathbf{P}\left( S_p-y\in [z, z+1], \min_{i\leq j\leq p} S_j\geq y, \tau_{\zeta_n}=i, \tau_{\zeta_n}^{(2)}> p, E_i^{(p)}\right)\nonumber\\
	& \lesssim p^a e^{-\lambda(mp)^b}\mathbf{P}\left(S_p-  \min_{i\leq j\leq p} S_j \leq z+1, E_i^{(p)}\right)\nonumber\\
	& \leq e^{-\alpha_p} \mathbf{P}\left(  \max_{0\leq j\leq p-i} S_j \leq z+1\right).
\end{align}
On the other hand, on the set $(E_i^{(p)})^c$, using \eqref{Step16}, 
\begin{align}\label{Step15}
	&  \mathbf{P}\left( \tau_{\zeta_n}=i, \tau_{\zeta_n}^{(2)}> p, (E_i^{(p)})^c \right) \nonumber\\
	& = \mathbf{P}\left(X<-\zeta_n\right) \mathbf{P}\left(|S_{p-1}-m(p-1)|> p^{(3-2b)/4}, \min_{j\leq p-1} X_j\geq -\zeta_n \right)\nonumber\\
	&\lesssim e^{-\alpha_n} n^{1-b+ \lambda b m^{b-1} A_1} \mathbf{P}\left(|S_{p-1}-m(p-1)|> p^{(3-2b)/4} \right). 
\end{align}
By Markov's inequality, under \eqref{Assumption4}, for any $Q\in \N$, there exists a constant $C(Q)\in(0,\infty)$ (for example, see \cite[p.60, Supplement 16]{Petrov1975}) such that for all $p\in \N$,
\begin{align}\label{Step17}
	& \mathbf{P}\left(|S_{p-1}-m(p-1)|> p^{(3-2b)/4} \right)\leq \frac{1}{p^{(3-2b)Q/4}} \mathbf{E}\left(\left| S_{p-1}-m(p-1)\right| ^Q \right)\nonumber\\
	& \leq \frac{C(Q)}{p^{(3-2b)Q/4}} (p-1)^{Q/2}.
\end{align}
Noticing that  $(3-2b)/4 > 1/2$ when $b\in (0, \tfrac{1}{2})$, fixing any $Q\in \N$ such that $(3-2b)Q/4 -Q/2 > 2(3-b+\lambda b m^{b-1}A_1)$, it follows from \eqref{Step15},  \eqref{Step17} and inequality $p^2\geq  n$  that 
\begin{align}\label{Step18}
	&  \mathbf{P}\left( \tau_{\zeta_n}=i, \tau_{\zeta_n}^{(2)}> p, (E_i^{(p)})^c \right)  \lesssim e^{-\alpha_n} n^{1-b+ \lambda b m^{b-1} A_1} p^{-2(3-b+\lambda b m^{b-1} A_1 ) }\lesssim \frac{1}{n^2 } e^{-\alpha_n}.
\end{align}
Therefore, (ii) follows directly from \eqref{Step14} and \eqref{Step18}.

\hfill$\Box$

\begin{lemma}\label{lem:local-limit-theorem}
	Assume \eqref{Assumption3b}  holds with $b\in (0,\tfrac{1}{2})$. Then 
	\begin{align}
		\sup_{ y<n^{(3-2b)/4} }\mathbf{P}(S_n-y\in[0,1])\lesssim ne^{-\alpha_n}. 
	\end{align}
\end{lemma}
\textbf{Proof: }  Noticing that we have the following upper bound:
\begin{align}\label{Step25}
	& \mathbf{P}(S_n-y\in[0,1])\nonumber\\
	& \leq 	\mathbf{P}(S_n\leq y+1, \tau_{\zeta_n}>n)+ 	\mathbf{P}\left(\tau_{\zeta_n}^{(2)}\leq n\right) + \sum_{i=1}^n 	\mathbf{P}(S_n-y\in[0,1], \tau_{\zeta_n}=i, \tau_{\zeta_n}^{(2)}>n)\nonumber\\
	&=: J_1+ J_2+ J_3.
\end{align}
We first treat $J_1$.  For any $y<n^{(3-2b)/4}$, combining Lemma \ref{lem1} (ii) and the inequality that $\theta_n n^{(3-2b)/4} \lesssim n^{b-1} n^{(3-2b)/4} =n^{(2b-1)/4}=o(1)$, we have
\begin{align}\label{Step24}
	J_1\lesssim e^{-\theta_n (mn-y-1)}\leq e^{\theta_n-\theta_n mn +\theta_n n^{(3-2b)/4}}\lesssim e^{-\theta_n mn}. 
\end{align}
According to the  definitions of $\theta_n$ and $\widehat{\zeta}_n$ in \eqref{Def-of-zeta}, 
\begin{align}
	e^{-\theta_n mn} & \lesssim e^{-\theta_n  m(n+1)} = \exp\left\{-\left(\lambda \widehat{\zeta}_n^{b-1}- \frac{A_2 \log \widehat{\zeta}_n}{\widehat{\zeta}_n}\right)\left(\widehat{\zeta}_n + A_1 n^{1-b}\log n\right) \right\}\nonumber\\
	& = \exp\left\{ -\lambda \widehat{\zeta}_n^b- \lambda A_1 \widehat{\zeta}_n^{b-1} n^{1-b}\log n + A_2 \log \widehat{\zeta}_n+o(1) \right\}\nonumber\\
	& = \exp\left\{-\lambda (mn)^b + \lambda b A_1 m^{b-1} \log n- \lambda A_1 \widehat{\zeta}_n^{b-1} n^{1-b}\log n + A_2 \log \widehat{\zeta}_n+o(1) \right\},
\end{align} 
where in the last equality we used the Taylor's expansion $\widehat{\zeta}_n^b= (mn)^b -   b A_1 m^{b-1} \log n+o(1)$. Noticing that 
\begin{align}
	& \lambda b A_1 m^{b-1} \log n- \lambda A_1 \widehat{\zeta}_n^{b-1} n^{1-b}\log n + A_2 \log \widehat{\zeta}_n\\
	& = -\left( (1-b)\lambda A_1 m^{b-1}-A_2+\varepsilon_n\right)\log n
\end{align}
for some $\varepsilon_n = o(1)$, 
according to our choice of $A_1$ in \eqref{Choice-of-A} and the definition of $\alpha_n$, when $n$ is large enough such that $a+ (1-b)\lambda A_1 m^{b-1}-A_2+\varepsilon_n \geq 3/2$, 
\begin{align}\label{Step20}
	e^{-\theta_n mn}& \lesssim e^{-\alpha_n} \exp\left\{-\left( a+ (1-b)\lambda A_1 m^{b-1}-A_2+\varepsilon_n \right)\log n \right\} \leq \frac{1}{n^{3/2}}e^{-\alpha_n}. 
\end{align}
Combining \eqref{Step24} and \eqref{Step20}, it holds that 
\begin{align}\label{Ineq-upper-of-J-1}
	J_1 \lesssim e^{-\alpha_n}. 
\end{align}
For $J_2$, by Lemma \ref{lem3}(i), $J_2$ is bounded from above by
\begin{align}\label{Ineq-upper-of-J-2}
	J_2 \lesssim n^2 \zeta_n^{2(a+1-b)} e^{-2\lambda\zeta_n^b} =e^{-\alpha_n} e^{-\lambda (mn)^b(1+o(1))}\lesssim e^{-\alpha_n}.
\end{align}
For $J_3$,  recall the defintion of $E_i^{(n)}$ in \eqref{Def-of-E-i-p}, then by \eqref{Step18} with $p=n$, we have
\begin{align}
	J_3 & \lesssim \frac{1}{n} e^{-\alpha_n}+ \sum_{i=1}^n 	\mathbf{P}\left(S_n-y\in[0,1], \tau_{\zeta_n}=i, \tau_{\zeta_n}^{(2)}>n, E_i^{(n)}\right)\nonumber\\
	& \leq e^{-\alpha_n}+n	\mathbf{E}\left(1_{E_n^{(n)}}\mathbf{P}(X+s-y\in [0,1] )|_{s=S_{n-1}} \right).
\end{align}
Noticing that on $E_n^{(n)}$, $y-S_{n-1}< 2n^{(3-2b)/4}- m(n-1)< -(mn-3n^{(3-2b)/4})$ for large $n$, $J_3$ is bounded from above by
\begin{align}\label{Ineq-upper-of-J-3}
	J_3 & \lesssim  e^{-\alpha_n}+n \sup_{z< -(mn-3n^{(3-2b)/4})} \mathbf{P}(X\in [z,z+1])\nonumber\\
	& \lesssim e^{-\alpha_n}+n \sup_{z< -(mn-3n^{(3-2b)/4})} \left(|z|^a e^{-\lambda |z|^b}\right)\asymp ne^{-\alpha_n}.
\end{align}
Therefore, combining \eqref{Step25}, \eqref{Ineq-upper-of-J-1}, \eqref{Ineq-upper-of-J-2} and \eqref{Ineq-upper-of-J-3}, we arrive at the desired result.

\hfill$\Box$

\begin{lemma}\label{lem:Renewal-theorem}
	Assume \eqref{Assumption3b} holds with $b\in (0,\tfrac{1}{2})$. Let $H(x,z)$ be a measurable function for $(x,z)\in [0,\infty) \times \R$ such that for any $\gamma>1$, there exists $K(\gamma)$ such that for all $x\geq 0$,
	\begin{align}
		\sup_{z\in \R} |H(x,z)| \leq K(\gamma) (1+x)^{-\gamma}.
	\end{align}
	Assume that the limit $H_\infty(x):= \lim_{z\to-\infty} H(x,z)$ exists for all $x\geq 0$. Then uniformly for any $|y|\leq 3n^{(3-2b)/4}$ and $|p-n|<n^{(3-2b)/4}$, 
	\begin{align}
		\lim_{n\to\infty} e^{\alpha_n} \mathbf{E}\left(1_{\{\tau_{\zeta_n}=p, S_p\geq y\}} H(S_p-y, X_p)\right)= \ell_\infty m^a \int_0^\infty H_\infty (x)\mathrm{d}x. 
	\end{align}
\end{lemma}
\textbf{Proof: } We only deal with the case that $H$ is non-negative since for the general case, we can decompose $H= \max\{H, 0\}- \max\{-H, 0\}$. Recall the definition of $E_i^{(p)}$ in \eqref{Def-of-E-i-p}. By \eqref{Step18}, it suffices to prove that 
\begin{align}\label{Step64}
	\lim_{n\to\infty} e^{\alpha_n} \mathbf{E}\left(1_{\{\tau_{\zeta_n}=p,   S_p\geq y\}} 1_{E_p^{(p)}} H(S_p-y, X_p)\right)=  \ell_\infty m^a \int_0^\infty H_\infty (z)\mathrm{d}z. 
\end{align}
Fixing any $\gamma>1$ such that $(3-2b)\gamma /4 > 1-b+\lambda b m^{b-1}A_1$, then by \eqref{Step16}, we have
\begin{align}
	& \limsup_{n\to\infty} e^{\alpha_n} \mathbf{E}\left(1_{\{\tau_{\zeta_n}=p,   S_p- y> p^{(3-2b)/4}\}} 1_{E_p^{(p)}} H(S_p-y, X_p)\right)\nonumber\\
	&\leq \limsup_{n\to\infty} e^{\alpha_n} \mathbf{P}(X<-\zeta_n) \frac{K(\gamma)}{(1+p^{(3-2b)/4})^\gamma}\nonumber\\
	&\lesssim \limsup_{n\to\infty} \frac{n^{1-b+\lambda b m^{b-1}A_1}}{n^{ (3-2b)\gamma /4}}=0. 
\end{align}
Therefore, to prove \eqref{Step64}, it remains to show that
\begin{align}\label{Step65}
	\lim_{n\to\infty} e^{\alpha_n} \mathbf{E}\left(1_{\{\tau_{\zeta_n}=p,   S_p-  y\in [0, p^{(3-2b)/4}]\}} 1_{E_p^{(p)}} H(S_p-y, X_p)\right)=  \ell_\infty m^a \int_0^\infty H_\infty (z)\mathrm{d}z. 
\end{align}
For any $\mathcal{A}\in \sigma(S_1,...,S_{p-1})$, according to the independence of $X_p$ and $\sigma(S_1,...,S_{p-1})$, it follows from \eqref{Assumption3b} that 
\begin{align}\label{Step66}
	& e^{\alpha_n} \mathbf{E}\left(1_{\mathcal{A}}1_{\{  S_p-  y\in [0, p^{(3-2b)/4}]\}} 1_{E_p^{(p)}} H(S_p-y, X_p)\right)\nonumber\\
	& =  \mathbf{E}\left( 1_{\mathcal{A}}1_{E_p^{(p)}} \int_{S_{p-1}-y +z\in [0,p^{(3-2b)/4}] } H(S_{p-1}+z-y, z)e^{\alpha_n} \ell(z) |z|^a e^{-\lambda |z|^b} \mathrm{d} z\right).
\end{align}
On $E_p^{(p)}\cap \left\{ S_{p-1}-y +z\in [0,p^{(3-2b)/4}] \right\}$, uniformly for all $|y|\leq 3n^{(3-2b)/4}$ and $|p-n|< n^{(3-2b)/4}$, when $n$ is large enough,
\[
|z+mn| \leq |S_{p-1}-m(p-1)|+ m(1+n^{(3-2b)/4})+ |y|+ p^{(3-2b)/4}\leq (8+m)n^{(3-2b)/4},
\]
which implies that uniformly on $E_p^{(p)}\cap \left\{ S_{p-1}-y +z\in [0,p^{(3-2b)/4}] \right\}$, $e^{\alpha_n}|z|^a e^{-\lambda |z|^b}= m^a (1+o(1))$. Plugging this back to \eqref{Step66}, we deduce that uniformly for all $|y|\leq 3n^{(3-2b)/4}$ and $|p-n|<n^{(3-2b)/4}$,  
\begin{align}\label{Step67}
	& e^{\alpha_n} \mathbf{E}\left(1_{\mathcal{A}}1_{\{  S_p-  y\in [0, p^{(3-2b)/4}]\}} 1_{E_p^{(p)}} H(S_p-y, X_p)\right)\nonumber\\
	& = \ell_\infty m^a (1+o(1)) \mathbf{E}\left( 1_{\mathcal{A}}1_{E_p^{(p)}} \int_{S_{p-1}-y +z\in [0,p^{(3-2b)/4}] } H(S_{p-1}+z-y, z)  \mathrm{d} z\right)\nonumber\\
	&= \ell_\infty m^a (1+o(1)) \mathbf{E}\left( 1_{\mathcal{A}}1_{E_p^{(p)}} \int_0^{p^{(3-2b)/4}} H(z, z+y-S_{p-1})  \mathrm{d} z\right).
\end{align}
Consequently, taking $\mathcal{A}:= \{\tau_{\zeta_n}\leq p-1\}$ in \eqref{Step67}, it follows from \eqref{Step16} that 
\begin{align}\label{Step68}
	& \limsup_{n\to\infty} e^{\alpha_n} \mathbf{E}\left(1_{\mathcal{A}}1_{\{  S_p-  y\in [0, p^{(3-2b)/4}]\}} 1_{E_p^{(p)}} H(S_p-y, X_p)\right)\nonumber\\
	&\lesssim  \limsup_{n\to\infty} \mathbf{P}(\tau_{\zeta_n}\leq p-1 ) \int_0^\infty \frac{K(\gamma)}{(1+z)^\gamma}\mathrm{d}z  \lesssim  \limsup_{n\to\infty}  n \mathbf{P}(X<-\zeta_n)=0. 
\end{align}
Noticing that when $n$ is large enough,  $\{\tau_{\zeta_n}\leq p\} \supset \left\{S_p-y\in [0, p^{(3-2b)/4}] \right\}\cap E_p^{(p)}$, combining \eqref{Step68} and \eqref{Step67}, we obtain that 
\begin{align}
	& \lim_{n\to\infty} e^{\alpha_n} \mathbf{E}\left(1_{\{\tau_{\zeta_n}=p,   S_p-  y\in [0, p^{(3-2b)/4}]\}} 1_{E_p^{(p)}} H(S_p-y, X_p)\right)\nonumber\\
	& = \lim_{n\to\infty} e^{\alpha_n} \mathbf{E}\left(1_{\{ S_p-  y\in [0, p^{(3-2b)/4}]\}} 1_{E_p^{(p)}} H(S_p-y, X_p)\right) \nonumber\\
	& = \ell_\infty m^a \lim_{n\to\infty}  \mathbf{E}\left( 1_{E_p^{(p)}} \int_0^{p^{(3-2b)/4}} H(z, z+y-S_{p-1})  \mathrm{d} z\right). 
\end{align}
Noticing that $z+y-S_{p-1}\leq -(mn -(4+m)n^{(3-2b)/4})$ and that $\mathbf{P}(E_p^{(p)})\to 1$, by dominated convergence theorem, we finally conclude that 
\begin{align}
	& \lim_{n\to\infty} e^{\alpha_n} \mathbf{E}\left(1_{\{\tau_{\zeta_n}=p,   S_p-  y\in [0, p^{(3-2b)/4}]\}} 1_{E_p^{(p)}} H(S_p-y, X_p)\right)= \ell_\infty m^a \int_0^\infty H_\infty (z)\mathrm{d}z,
\end{align}
which implies \eqref{Step65} and we arrive at the desired result.

\hfill$\Box$

In the rest part of this section, we gather some inequalities which will be used in Section \ref{S4} in the proof of Proposition \ref{prop:key-proposition}. 
For each $x\geq 0$, define
\begin{align}\label{Def-of-renewal-function}
	R(x):= \sum_{n=0}^\infty \mathbf{P}\left(\max_{0\leq j\leq n}S_j \leq x\right).
\end{align}
Since $\mathbf{E}(X)>0$, it is well-knwon that (for example, see \cite[Lemma 2.1, p.1950]{Gut2009}) $R(x)$ is finite.

\begin{lemma}\label{lem:bad-even5}
	For any $\varepsilon, L> 0$, there exists a constant $\Lambda=\Lambda(\varepsilon,L )\in (0,\infty)$ such that for all $y <n^{(3-2b)/4}$ and $z\in [0,L]$, when $n$ is large enough, 
	\begin{align}
		\mathbf{P}\left(S_n-y\in [z,z+1], \min_{\tau_{\zeta_n}\leq j\leq n}S_j\geq y, |S_n-S_{\tau_{\zeta_n}}| >\Lambda, \tau_{\zeta_n} \leq n\right)\leq \varepsilon e^{-\alpha_n}.
	\end{align}
\end{lemma}
\textbf{Proof: } Recall the definition of $E_i^{(n)}$ in \eqref{Def-of-E-i-p}. Combining \eqref{Step11} and \eqref{Step12}, for any $1\leq i\leq n$,
\begin{align}\label{Step31}
	&\mathbf{P}\left(|S_n-S_i|>\Lambda, S_n-y\in [z, z+1], E_i^{(n)}\right)\nonumber\\
	&= \mathbf{E}\left(1_{\{|S_n-S_i|>\Lambda\}}1_{E_i^{(n)}} \mathbf{P}\left(X_i+t-y \in [z,z+1]\right)\bigg|_{t=S_n-X_i}\right)\nonumber\\
	&\lesssim  \mathbf{P}\left(|S_n-S_i|>\Lambda, E_i^{(n)}\right) n^a e^{-\lambda(mn)^b} \leq e^{-\alpha_n}  \mathbf{P}\left(|S_n-S_i|>\Lambda\right) .
\end{align}
Also, from Lemma \ref{lem3} (ii), for each $1\leq T\leq n$, 
\begin{align}\label{Step32}
	& \mathbf{P}\left(S_n-y\in [z,z+1], \min_{\tau_{\zeta_n}\leq j\leq n} S_j\geq y, \tau_{\zeta_n} \leq n-T, \tau_{\zeta_n}^{(2)}>n \right)\nonumber\\
	& \lesssim \sum_{i=1}^{n-T} \left( \mathbf{P}\left(\max_{0\leq j\leq n-i} S_j \leq z+1\right) +\frac{1}{n^2}\right)e^{-\alpha_n}\nonumber\\
	&\leq e^{-\alpha_n}\left(\sum_{i=T}^\infty \mathbf{P}\left(\max_{0\leq j\leq i} S_j \leq z+1\right) +\frac{1}{n} \right). 
\end{align}
Combining Lemma \ref{lem3}(i),  \eqref{Step18}, \eqref{Step31} and \eqref{Step32}, we have
\begin{align}\label{Step33}
	&\mathbf{P}\left(S_n-y\in [z,z+1], \min_{\tau_{\zeta_n}\leq j\leq n}S_j\geq y, |S_n-S_{\tau_{\zeta_n}}| >\Lambda, \tau_{\zeta_n} \leq n\right)\nonumber\\
	&\lesssim n^2  \zeta_n^{2(a+1-b)} e^{-2\lambda\zeta_n^b} + \sum_{i=1}^n \frac{1}{n^2 } e^{-\alpha_n} + e^{-\alpha_n}\left(\sum_{i=T}^\infty \mathbf{P}\left(\max_{0\leq j\leq i} S_j \leq z+1\right) +\frac{1}{n} \right)\nonumber\\
	&\quad + \sum_{i=n-T+1}^n \mathbf{P}\left(|S_n-S_i|>\Lambda, S_n-y\in [z, z+1], E_i^{(n)}\right)\nonumber\\
	&\lesssim e^{-\alpha_n}\left(\sum_{i=T}^\infty \mathbf{P}\left(\max_{0\leq j\leq i} S_j \leq z+1\right) +\frac{1}{n} \right) + e^{-\alpha_n}\sum_{i=1}^{T-1} \mathbf{P}\left(|S_i|>\Lambda\right). 
\end{align} 
Therefore, there exists a constant $C$ such that 
\begin{align}
	&\mathbf{P}\left(S_n-y\in [z,z+1], \min_{\tau_{\zeta_n}\leq j\leq n}S_j\geq y, |S_n-S_{\tau_{\zeta_n}}| >\lambda, \tau_{\zeta_n} \leq n\right)\nonumber\\
	&\leq e^{-\alpha_n} C\left(\sum_{i=T+1}^\infty \mathbf{P}\left(\max_{0\leq j\leq i} S_j \leq L+1\right) +\frac{1}{n} + \sum_{i=1}^{T-1} \mathbf{P}\left(|S_i|>\lambda\right) \right).
\end{align}
By \eqref{Def-of-renewal-function}, for any $\varepsilon, L>0$, let $T$ and $\lambda$ be large enough such that $\sum_{i=T+1}^\infty \mathbf{P}\left(\max_{0\leq j\leq i} S_j \leq L+1\right)< \varepsilon / (3C)$ and $\sum_{i=1}^{T-1} \mathbf{P}\left(|S_i|>\lambda\right) <\varepsilon/(3C)$, then we complete the proof of the lemma by taking $n> \max\{T, (3C)/\varepsilon\}$.

\hfill$\Box$

\begin{lemma}\label{lem4}
	For each $L>0$ and $T\in \N$, when $n$ is large enough, for any $n^{(3-2b)/4}\leq k\leq n$, it holds that 
	\begin{align}
		& \sup_{|z+mk|\leq 3n^{(3-2b)/4}} \mathbf{P}\left(S_{n-k}\leq z, \min_{\tau_{\zeta_n}\leq i\leq n-k} S_i\geq z-L, \tau_{\zeta_n}\in [n-k-T, n-k], \tau_{\zeta_n}^{(2)}>n-k\right)\nonumber\\
		&\lesssim_{L,T}e^{-\alpha_n}. 
	\end{align} 
\end{lemma}
\textbf{Proof: }  Set $k_0:= \max\{n-k-T, 1\}$. Firstly noticing that 
\begin{align}
	& \mathbf{P}\left(S_{n-k}\leq z, \min_{\tau_{\zeta_n}\leq i\leq n-k} S_i\geq z-L, \tau_{\zeta_n}\in [n-k-T, n-k], \tau_{\zeta_n}^{(2)}>n-k\right)\nonumber\\
	&\leq \sum_{j=k_0}^{n-k} \mathbf{P}\left(S_{n-k}\leq z, \min_{j\leq i\leq n-k} S_i\geq z-L, \tau_{\zeta_n}=j,  \tau_{\zeta_n}^{(2)}>n-k\right)\nonumber\\
	&\leq \sum_{j=k_0}^{n-k} \mathbf{P} \left(X_j <-\zeta_n\right) \mathbf{P}\left( |S_{n-k}-X_j-m(n-k-1)|>n^{(3-2b)/4} \right)\nonumber\\
	&\quad + \sum_{j=k_0}^{n-k} \mathbf{E}\left(1_{\{|S_{n-k}-X_j-m(n-k-1)|\leq n^{(3-2b)/4} \}}\mathbf{P}(X+t\in [z-L,z])\bigg|_{t=S_{n-k}-X_j}\right).
\end{align}
Since on $\{|S_{n-k}-X_j-m(n-k-1)|\leq n^{(3-2b)/4} \}$, for $t=S_{n-k}-X_j$ and $|z+mk|\leq 3n^{(3-2b)/4}$, we have
$z-t\leq 3n^{(3-2b)/4}-mk-m(n-k-1)+n^{(3-2b)/4}\leq -(mn-5n^{(3-2b)/4})$ when $n$ is large enough. Therefore, together with \eqref{Step16}, the above inequality is bounded by
\begin{align}\label{Step34}
	& \mathbf{P}\left(S_{n-k}\leq z, \min_{\tau_{\zeta_n}\leq i\leq n-k} S_i\geq z-L, \tau_{\zeta_n}\in [n-k-T, n-k], \tau_{\zeta_n}^{(2)}>n-k\right)\nonumber\\
	&\lesssim_L (T+1)e^{-\alpha_n} n^{1-b+ \lambda b m^{b-1} A_1} \mathbf{P} \left(|S_{n-k-1}-m(n-k-1)|>n^{(3-2b)/4}  \right)\nonumber\\
	&\quad + (T+1) \sup_{z< - (mn-5n^{(3-2b)/4})} |z|^a e^{-\lambda |z|^b}.
\end{align}
Let $Q$ be the fixed integer in \eqref{Step17}, then combining  \eqref{Step17} and \eqref{Step34}, we conclude that 
\begin{align}
	& \mathbf{P}\left(S_{n-k}\leq z, \min_{\tau_{\zeta_n}\leq i\leq n-k} S_i\geq z-L, \tau_{\zeta_n}\in [n-k-T, n-k], \tau_{\zeta_n}^{(2)}>n-k\right)\nonumber\\
	&\lesssim_{L, T}  e^{-\alpha_n} n^{1-b+ \lambda b m^{b-1} A_1} \frac{(n-k-1)^{Q/2}}{n^{(3-2b)Q/4}}+ n^ae^{-\lambda (mn)^b} \nonumber\\
	&  \leq e^{-\alpha_n} n^{1-b+ \lambda b m^{b-1} A_1} \frac{n^{Q/2}}{n^{(3-2b)Q/4}}+ e^{-\alpha_n}\lesssim e^{-\alpha_n}. 
\end{align}
We are done.

\hfill$\Box$

\section{Proofs of the Theorems \ref{thm1} and \ref{thm2}}\label{S3}

This section is devoted to proving Theorems \ref{thm1} and \ref{thm2}. 

Recall that $\{W_n = \sum_{|u|=n}e^{-V(u)}, \mathcal{F}_n, n\in \N,  \P\}$ is the additive martingale, where $\mathcal{F}_n$ is the natural filtration of the branching random walk up to generation $n$. Define 
\begin{align}\label{Change-of-measure}
	\frac{\mathrm{d} \Q }{\mathrm{d} \P}\bigg|_{\mathcal{F}_n}:= W_n. 
\end{align}
Denote by $\widehat{\mathscr{L}}$ the  law of $\mathscr{L}$ under $\Q$. Lyons \cite{Lyons} gave the following description of
the law of the branching random walk under $\Q$: there is a spine process denoted by  $\{w_n\}_{n\ge 0}$ with $w_0 =\emptyset$ and the initial position of the spine  is $V(w_0)=0$. At generation $n=1$, $w_0$ dies and splits into a random number of offspring equal in law to $\widehat{\mathscr{L}}$.  Choose one offspring $x$ from all the offspring of $w_0$ with probability proportional to $e^{-V(x)}$, and call it $w_1$. $w_1$ evolves independently as $w_0$ and the other unmarked offspring evolve  independently as in the original branching random walk. By Lyons \cite{Lyons}, for any $u\in \mathbb{T}$ with $|u|=n$, we have
\begin{align}\label{Gibbs-measure}
	\Q\left(w_n =u |\mathcal{F}_n\right) = \frac{e^{-V(u)}}{W_n}.
\end{align}
Moreover,  the position process $\left\{V(w_n)\right\}_{n\geq 0}$ along the spine under $\Q$ is equal in the law to $\left\{S_n\right\}_{n\geq 0}$ defined in \eqref{Law-of-spine}. According to \eqref{Gibbs-measure}, we have the following many-to-one formula: for any Borel measurable function $f$, it holds that 
\begin{align}\label{Many-to-one}
	\mathbf{E}\left(f\left(S_1,...,S_n\right)\right)=\E\Big(\sum_{|u|=n} f\left(V(u_1),...,V(u)\right)e^{-V(u)}\Big).
\end{align}
For each $u\in \mathbb{T}$, define 
\begin{align}
	\tau_{\zeta_n}^{(u)}& := \inf\left\{1\leq i \leq |u|: V(u_i)- V(u_{i-1})<-\zeta_n \right\},\nonumber\\
	\tau_{\zeta_n}^{(2, u)}& := \inf\left\{ 	\tau_{\zeta_n}^{(u)}< i \leq |u|: V(u_i)- V(u_{i-1})<-\zeta_n \right\},
\end{align}
with the convention $\inf \emptyset := \infty$.  To simplify the notation, set
\begin{align}\label{Def-of-T-n}
	\mathbb{T}_n:= \left\{u\in \mathbb{T}: |u|=n\right\}.
\end{align} 

\begin{lemma}\label{lem:bad-event1}
	For any $\varepsilon>0$, there exists an integer $n_0(\varepsilon)$ such that for all $n\geq n_0(\varepsilon)$ and $x\geq 0$, 
	\begin{align}\label{Eq-at-least-one-jump}
		\P\left(\exists u\in \mathbb{T}_n, V(u)\leq \alpha_n-x, \tau_{\zeta_n}^{(u)}>n\right)\leq \varepsilon e^{-x}
	\end{align}
	and 
	\begin{align}\label{Eq-at-most-two-jumps-2}
		\P\left(\exists u\in \mathbb{T}_n, V(u)\leq \alpha_n-x,  \tau_{\zeta_n}^{(2, u)}\leq n\right)\leq \varepsilon e^{-x}.
	\end{align}
\end{lemma}
\textbf{Proof: } Combining the union bound and many-to-one formula \eqref{Many-to-one}, 
\begin{align}
	& \P\left(\exists u\in \mathbb{T}_n, V(u)\leq \alpha_n-x, \tau_{\zeta_n}^{(u)}>n\right)\nonumber\\
	& \leq \E\Big( \sum_{|u|=n} 1_{\{V(u)\leq \alpha_n-x, \tau_{\zeta_n}^{(u)}>n \}} \Big) =\mathbf{E}\left(e^{S_n}1_{\{ S_n\leq \alpha_n-x, \tau_{\zeta_n}>n \}}\right)\nonumber\\
	&\leq e^{\alpha_n-x} \mathbf{P}\left( S_n\leq \alpha_n-x, \tau_{\zeta_n}>n \right).
\end{align}
By Lemma \ref{lem1} (ii),  the above probability is bounded from above by
\begin{align}\label{Step19}
	& \P\left(\exists u\in \mathbb{T}_n, V(u)\leq \alpha_n-x, \tau_{\zeta_n}^{(u)}>n\right)\lesssim e^{\alpha_n-x}e^{-\theta_n \left( mn+x-\alpha_n \right)} \leq e^{\alpha_n-x} e^{-\theta_n \left( mn-\alpha_n \right)},\quad 
\end{align}
where $\theta_n$ is defined in \eqref{Def-of-zeta}.  Since $\theta_n \alpha_n\asymp n^{b-1} n^b=o(1)$, combining  \eqref{Step20} and \eqref{Step19}, we get
\begin{align}
	& \P\left(\exists u\in \mathbb{T}_n, V(u)\leq \alpha_n-x, \tau_{\zeta_n}^{(u)}>n\right)\lesssim e^{\alpha_n-x} e^{-\theta_n mn }\lesssim \frac{1}{n^{3/2}}e^{-x},
\end{align}
which implies \eqref{Eq-at-least-one-jump}.

For \eqref{Eq-at-most-two-jumps-2}, combining the union bound and many-to-one formula \eqref{Many-to-one},
\begin{align}
	& \P\left(\exists u\in \mathbb{T}_n, V(u)\leq \alpha_n-x,  \tau_{\zeta_n}^{(2, u)}\leq n\right) \leq \E\Big(\sum_{|u|=n} 1_{\left\{V(u)\leq \alpha_n-x,  \tau_{\zeta_n}^{(2, u)}\leq n\right\}}\Big)\nonumber\\
	& = \mathbf{E}\Big(e^{S_n} 1_{\left\{S_n \leq \alpha_n-x,  \tau_{\zeta_n}^{(2)}\leq n\right\}}\Big)\nonumber\\
	& \leq  e^{\alpha_n-x} \mathbf{P}\left( \tau_{\zeta_n}^{(2)}\leq n\right).
\end{align}
Combining Lemma \ref{lem3} (i) and the above inequality, we deduce that
\begin{align}
	& \P\left(\exists u\in \mathbb{T}_n, V(u)\leq \alpha_n-x, \tau_{\zeta_n}^{(2, u)}\leq n\right)\nonumber\\
	&\lesssim e^{\alpha_n-x} n^2 \zeta_n^{2(a+1-b)} e^{-2\lambda\zeta_n^b}  = e^{-x}e^{-\lambda (mn)^b(1+o(1))}, 
\end{align} 
which implies \eqref{Eq-at-most-two-jumps-2}.

\hfill$\Box$

To simplify the notation, for each $x, L$ and $n$, set
\begin{align}\label{Def-of-y}
	y_n= y_n(x,L):=\alpha_n-x-L\quad \mbox{and}\quad \mathcal{H}_{\zeta_n, p}^u:= \left\{ \tau_{\zeta_n}^{(u)}\leq p < \tau_{\zeta_n}^{(2, u)}\right\}.
\end{align}

\begin{lemma}\label{lem:bad-event2}
	For any $L_0\in \N$, when $n$ is large enough, for any $x\geq 0$, and any $L\in [L_0, n^{(2b+1)/4}]$, it holds that 
	\begin{align}\label{Eq1}
		& \P\bigg(\exists u\in \mathbb{T}_n, V(u)\leq \alpha_n -x, \min_{\tau_{\zeta_n}^{(u)}\leq j\leq n} V(u_j)-y_n\in [0, 1], \mathcal{H}_{\zeta_n, n}^u\bigg)\nonumber\\
		&\lesssim (1+L)^{-2}e^{-x}.
	\end{align}
	As a consequence,  for any $L_0\in \N$, when $n$ is large enough, for any $L\geq L_0$, we have
	\begin{align}\label{Eq2}
		& \P\bigg(\exists u\in \mathbb{T}, |u|=n, V(u)\leq \alpha_n -x, \min_{\tau_{\zeta_n}^{(u)}\leq j\leq n} V(u_j)\leq y_n, \mathcal{H}_{\zeta_n, n}^u \bigg)\nonumber\\
		&\lesssim (1+L)^{-1}e^{-x}.
	\end{align}
\end{lemma}
\textbf{Proof: } 
On the set $\min_{\tau_{\zeta_n}^{(u)}\leq j\leq n} V(u_j)-y_n\in [0, 1]$, there exists $p\in [\tau_{\zeta_n}^{(u)}, n]$ such that $V(u_p)-y_n\in [0,1]$ and that $\min_{\tau_{\zeta_n}^{(u)}\leq j\leq p} V(u_j)\geq y_n$. Therefore,  $\tau_{\zeta_n}^{(u)}=\tau_{\zeta_n}^{(u_p)}$ and this together with $\alpha_n -x= y_n+L$ implies that 
\begin{align}\label{Step27}
	& \P\bigg(\exists u\in \mathbb{T}_n, V(u) \leq \alpha_n-x, \min_{\tau_{\zeta_n}^{(u)}\leq j\leq n} V(u_j)-y_n\in [0, 1], \mathcal{H}_{\zeta_n, n}^u \bigg)\nonumber\\
	& \leq \sum_{p=1}^n \P\bigg( \exists u\in \mathbb{T}_n, V(u)-y_n\leq L, \min_{\tau_{\zeta_n}^{(u)}\leq j\leq n} V(u_j)\geq y_n, V(u_p)-y_n\in [0,1],  \tau_{\zeta_n}^{(u_p)}\leq p < \tau_{\zeta_n}^{(2, u)}
	\bigg)\nonumber\\
	&\leq \sum_{p=1}^{n-2[L^{4/(3-2b)}]} I_1(p)+ \sum_{p=n-2[L^{4/(3-2b)}]+1}^{n} I_2(p),
\end{align}
where 
\begin{align}
	I_1(p)& :=  \P\bigg(\exists u\in \mathbb{T}_n, V(u) -y_n\leq L, \min_{\tau_{\zeta_n}^{(u)}\leq j\leq n} V(u_j)\geq y_n, V(u_p)-y_n\in[0,1],  \mathcal{H}_{\zeta_n, p}^u \bigg),\nonumber\\
	I_2(p) &:= \P\bigg(\exists u\in \mathbb{T}_p, V(u) -y_n\in [0,1], \min_{\tau_{\zeta_n}^{(u)}\leq j\leq p} V(u_j)\geq y_n,  \mathcal{H}_{\zeta_n, p}^u \bigg).
\end{align}
We treat $I_2(p)$ first. For $n-2[L^{4/(3-2b)}]+1\leq p\leq n$, according to the union bound and many-to-one formula \eqref{Many-to-one}, we have
\begin{align}\label{Step22}
	I_2(p) &\leq  \E\bigg(\sum_{|u|=p} 1_{\left\{V(u) -y_n\in [0,1], \min_{\tau_{\zeta_n}^{(u)}\leq j\leq p} V(u_j)\geq y_n,  \mathcal{H}_{\zeta_n, p}^u\right\}} \bigg)\nonumber\\
	& = \mathbf{E}\left(e^{S_p}1_{\left\{ S_p -y_n\in [0,1], \min_{\tau_{\zeta_n}\leq j\leq p} S_j\geq y_n,  \tau_{\zeta_n}\leq p < \tau_{\zeta_n}^{(2)}  \right\}}\right)\nonumber\\
	&\leq e^{y_n+1}\mathbf{P}\left( S_p -y_n\in [0,1], \min_{\tau_{\zeta_n}\leq j\leq p} S_j\geq y_n,  \tau_{\zeta_n}\leq p < \tau_{\zeta_n}^{(2)}  \right). 
\end{align} 
Since $2L^{4/(3-2b)}\leq 2n^{(2b+1)/(3-2b)}< n/2$ for large $n$, we have $p\geq n/2 >\sqrt{n}$. Also, noticing that $y_n\leq \alpha_n<n^{(2b+1)/4} < 4p^{(3-2b)/4}$ for large $n$ since $b< (2b+1)/4 <1/2< (3-2b)/4$, by Lemma \ref{lem3} (ii), we conclude from the above inequality 
\begin{align}\label{Ineq-Upper-I-2-p'}
	\sum_{p=n-2[L^{4/(3-2b)}]+1}^{n} I_2(p) & \lesssim e^{y_n}  \sum_{p=n-2[L^{4/(3-2b)}]+1}^{n} \sum_{i=1}^p \left( \mathbf{P}\left(\max_{0\leq j\leq p-i} S_j \leq 1\right) e^{-\alpha_p} +\frac{1}{n^2}e^{-\alpha_n} \right) \nonumber\\
	& \lesssim e^{y_n} \sum_{p=n-2[L^{4/(3-2b)}]+1}^{n} \left(e^{-\alpha_p}+e^{-\alpha_n}\right)\nonumber\\
	&= e^{-x-L} \sum_{p=n-2[L^{4/(3-2b)}]+1}^{n} \left(e^{\alpha_n-\alpha_p}+1\right).
\end{align}
Since $0\leq \log n- \log p \leq \log n - \log (n-2n^{(2b+1)/(3-2b)} +1) =  o(1) $ for all $n-2[L^{4/(3-2b)}]+1\leq p\leq n$ and $1\leq L\leq  n^{(2b+1)/4}$, according to Taylor's expansion and the fact that $n^{b-1}L^{(1+2b)/(3-2b)}\leq n^{b-1}L\leq n^{b-1}n^{(2b+1)/4}= n^{3(2b-1)/4}= o(1)$, we conclude that that for large $n$,
\[
e^{\alpha_n -\alpha_p} \lesssim e^{\lambda m^b (n^b -p^b)} \leq e^{\lambda m^b (n^b -(n-2L^{4/(3-2b)})^b)} = e^{O(n^{b-1}L^{4/(3-2b)})} = e^{o(L)} \leq e^{L/2}. 
\]
Therefore, it follows from \eqref{Ineq-Upper-I-2-p'} that 
\begin{align}\label{Ineq-Upper-I-2-p}
	&\sum_{p=n-2[L^{4/(3-2b)}]+1}^{n} I_2(p)  \lesssim e^{-x-L} \sum_{p=n-2[L^{4/(3-2b)}]+1}^{n} e^{L/2}\leq 2(L+1)^{4/(3-2b)}e^{-x-L/2}.
\end{align}
Now we treat $I_1(p)$. Similarly combining the union bound and many-to-one formula \eqref{Many-to-one}, 
\begin{align}
	I_1(p) & \leq \E\bigg(\sum_{|u|=n} 1_{\left\{V(u) -y_n\leq L, \min_{\tau_{\zeta_n}^{(u)}\leq j\leq n} V(u_j)\geq y_n, V(u_p)-y_n\in[0,1],  \mathcal{H}_{\zeta_n, p}^u \right\}} \bigg)\nonumber\\
	& = \mathbf{E}\left(e^{S_n}1_{\left\{ S_n -y_n\leq L, \min_{\tau_{\zeta_n} \leq j\leq n} S_j \geq y_n, S_p-y_n\in[0,1],  \tau_{\zeta_n} \leq p < \tau_{\zeta_n}^{(2)}  \right\}}\right)\nonumber\\
	& \leq e^{y_n+L}\mathbf{P} \Big(S_n -y_n\leq L, \min_{\tau_{\zeta_n} \leq j\leq n} S_j \geq y_n, S_p-y_n\in[0,1],  \tau_{\zeta_n} \leq p < \tau_{\zeta_n}^{(2)}  \Big). 
\end{align}
Applying the Markov property at time $p$ on the right hand side of the above probability, we get
\begin{align}\label{upp-of-I-1}
	I_1(p)  & \leq  e^{\alpha_n -x} \mathbf{E}\bigg(1_{\left\{ \min_{\tau_{\zeta_n} \leq j\leq p} S_j \geq y_n, S_p-y_n\in[0,1],  \tau_{\zeta_n} \leq p < \tau_{\zeta_n}^{(2)} \right\}}\nonumber\\
	&\quad\quad\times \mathbf{P}\left( \min_{0\leq j\leq n-p} S_j\geq y_n -z, S_{n-p}-y_n\leq L-z, \tau_{\zeta_n}> n-p\right)\bigg|_{z=S_p}\bigg)\nonumber\\
	& \leq e^{\alpha_n -x} \mathbf{P}\left( \min_{\tau_{\zeta_n} \leq j\leq p} S_j \geq y_n, S_p-y_n\in[0,1],  \tau_{\zeta_n} \leq p < \tau_{\zeta_n}^{(2)}\right) \nonumber\\
	&\quad \times \mathbf{P}\left( \min_{0\leq j\leq n-p} S_j\geq -1,  S_{n-p}\in [-1, L+1], \tau_{\zeta_n}> n-p\right).
\end{align}
When $p\leq [n^{(3-2b)/4}]$, applying Lemma \ref{lem1} (ii) with $q=n-p$, we have
\begin{align}
	I_1(p) &\leq e^{\alpha_n -x} \mathbf{P}\left(   S_{n-p} - m(n-p)\leq  L+1- m(n-p), \tau_{\zeta_n}>n-p \right)\nonumber\\
	&\lesssim e^{\alpha_n-x}e^{-\theta_n (m(n-p)-L-1)}\lesssim e^{\alpha_n-x} e^{-\theta_n mn},
\end{align}
where in the last inequality we used the fact that $\theta_n (p+L)\lesssim n^{b-1}(n^{(3-2b)/4}+ n^{(2b+1)/4})=o(1)$. Therefore, combining \eqref{Step20} and inequality $(L+1)^2 n^{(3-2b)/4} \lesssim n^{(2b+1)/2}n^{(3-2b)/4} = n^{3/2- (1-2b)/4}\leq n^{3/2}$, we obtain
\begin{align}\label{Ineq-Upper-I-1-p-small}
	\sum_{p=1}^{[n^{(3-2b)/4}]} I_1(p) \lesssim \frac{n^{(3-2b)/4}}{n^{3/2}}e^{-x}\lesssim (1+L)^{-2}e^{-x}.
\end{align}

When $[n^{(3-2b)/4}]< p \leq n-2[L^{4/(3-2b)}]$, we have $p\in [\sqrt{n}, n]$ and $L< (n-p)^{(3-2b)/4}$. Therefore, combining Lemma \ref{lem3} (ii) and Lemma \ref{lem:local-limit-theorem}, 
\begin{align}\label{Step26}
	I_1(p) & \lesssim e^{\alpha_n-x}\sum_{i=1}^p \left(\mathbf{P}\left(\max_{j\leq p-i+1} S_j \leq 1\right) e^{-\alpha_p} +\frac{1}{n^2}e^{-\alpha_n} \right)\times \mathbf{P}(S_{n-p}\in [-1, L+1])\nonumber\\
	& \lesssim e^{\alpha_n -x} (1+L) \left(\sum_{i=1}^\infty\mathbf{P}\left(\max_{j\leq i} S_j \leq 1\right) e^{-\alpha_p} + \frac{1}{n}e^{-\alpha_n}  \right)  (n-p)e^{-\alpha_{n-p}} \nonumber\\
	& \lesssim e^{\alpha_n -x} (1+L)e^{-\alpha_p} (n-p)e^{-\alpha_{n-p}} + e^{-x} (1+L) e^{-\alpha_{n-p}} =: I_{11}(p)+ I_{12}(p),
\end{align}
where in the second inequality we used the fact that the function $R$ defined in \eqref{Def-of-renewal-function} is finite.  
For $I_{12}$, it is easy to see that 
\begin{align}\label{upp-of-I-12}
	\sum_{p=1}^{n-2[L^{4/(3-2b)}]} I_{12}(p) & = e^{-x}(1+L) \sum_{p=2[L^{4/(3-2b)}]}^{n-1} e^{-\alpha_p} \nonumber\\
	& \leq e^{-x}(1+L) \sum_{p=2[L^{4/(3-2b)}]}^{\infty} p^{|a|} e^{- \lambda (mp)^b} \lesssim e^{-x} (1+L)^{-2}. 
\end{align}
If $p\leq n- [n^{(3-2b)/4}]-1 $, then $n^{(3-2b)/4}<  p < n- n^{(3-2b)/4}$. In this case, $p^b +(n-p)^b$ is increasing when $n^{(3-2b)/4}<p <n/2$ and decreasing when $n/2<p < n- n^{(3-2b)/4}$, which implies  that 
\begin{align}
	& \sum_{p= [n^{(3-2b)/4}]+1}^{ (n- [n^{(3-2b)/4}]-1 )\land (n- 2[L^{4/(3-2b)}])} I_{11}(p)\nonumber\\
	& \lesssim  n^{3|a|+1}(L+1)e^{-x}\sum_{p= [n^{(3-2b)/4}]+1}^{n- [n^{(3-2b)/4}]-1 } e^{ \lambda(mn)^b }  e^{-\lambda (m n^{(3-2b)/4})^b - \lambda m^b(n- n^{(3-2b)/4})^b} \nonumber\\
	&= e^{-x}(L+1)n^{3|a|+2} e^{\lambda (mn)^b}   e^{-\lambda (m n^{(3-2b)/4})^b - \lambda  m^b(n- n^{(3-2b)/4})^b}.
\end{align}
According to Taylor's expansion, $\lambda (mn)^b - \lambda  m^b(n- n^{(3-2b)/4})^b = \lambda (mn)^{b-1} n^{(3-2b)/4} + o(1)= o(1)$. Therefore, we get that 
\begin{align}\label{Ineq-Upper-I-1-p-middle}
	\sum_{p= [n^{(3-2b)/4}]+1}^{ (n- [n^{(3-2b)/4}]-1 )\land (n- 2[L^{4/(3-2b)}])} I_{11}(p)& \lesssim e^{-x}(L+1) n^{3|a|+2}  e^{-\lambda (m n^{(3-2b)/4})^b } \lesssim (L+1)^{-2}e^{-x}. \quad 
\end{align}
Now we treat the last case $n- [n^{(3-2b)/4}]-1 < p \leq n- 2[L^{4/(3-2b)}]$. In this case,  uniformly, 
\begin{align}
	\alpha_n -\alpha_p =-a\log (n/p)+ \lambda (mn)^b -\lambda (mn-p)^b = o(1),
\end{align}
which together with the definition of $I_{11}(p)$ in \eqref{Step26} implies that 
\begin{align}\label{Ineq-Upper-I-1-p-large}
	\sum_{p=n-[n^{(3-2b)/4}]}^{n-2[L^{4/(3-2b)}]}I_{11}(p)&\lesssim  (L+1)e^{-x}\sum_{p=n-[n^{(3-2b)/4}]}^{n-2[L^{4/(3-2b)}]} e^{\alpha_n } e^{-\alpha_p} (n-p)e^{-\alpha_{n-p}}\nonumber\\&\lesssim e^{-x}(L+1)\sum_{p=n-[n^{(3-2b)/4}]}^{n-2[L^{4/(3-2b)}]}   (n-p)e^{-\alpha_{n-p}} = e^{-x}(L+1)\sum_{p=2[L^{4/(3-2b)}]}^{[n^{(3-2b)/4}]}   p e^{-\alpha_{p}} \nonumber\\
	& \leq e^{-x}(L+1)\sum_{p=2[L^{4/(3-2b)}]}^{\infty}   p^{|a|+1} e^{-\lambda (mp)^b} \lesssim (L+1)^{-2}e^{-x}.
\end{align}
Combining  \eqref{Ineq-Upper-I-1-p-small}, \eqref{upp-of-I-12}, \eqref{Ineq-Upper-I-1-p-middle} and \eqref{Ineq-Upper-I-1-p-large}, we get that $\sum_{p=1}^{n-2[L^{4/(3-2b)}]} I_1(p)\lesssim (1+L)^{-2} e^{-x}$. Combining  \eqref{Step27}, \eqref{Ineq-Upper-I-2-p} and the above inequality, we get \eqref{Eq1}. 

Now we are going to prove \eqref{Eq2}. 
Noticing that by \eqref{Many-to-one}, for any $z\geq 0$,
\begin{align}\label{Tail-lowest}
	& \P\left(\exists u\in \mathbb{T}: V(u)< -z\right)\leq \sum_{k=1}^\infty \P\left(\exists u\in \mathbb{T}_k: V(u)< -z, \min_{j\leq k} V(u_j)\geq -z\right)\nonumber\\
	& \leq \sum_{k=1}^\infty \E\bigg(\sum_{|u|=k} 1_{\left\{ V(u)< -z, \min_{j\leq k} V(u_j)\geq -z \right\}} \bigg) = \sum_{k=1}^\infty \mathbf{E} \left(e^{S_k} 1_{\left\{ S_k< -z, \min_{j\leq k} S_j\geq -z \right\}} \right) \nonumber\\
	&\leq e^{-z}  \sum_{k=1}^\infty \mathbf{P} \Big( S_k< -z, \min_{j\leq k} S_j\geq -z \Big) =e^{-z}.
\end{align}
If $L\geq 2\alpha_n$, then $\alpha_n -L\leq -L/2$ and  by  \eqref{Tail-lowest}, we have
\begin{align}\label{Step28}
	& \P\bigg(\exists u\in \mathbb{T}_n, V(u)\leq \alpha_n -x, \min_{\tau_{\zeta_n}^{(u)}\leq j\leq n} V(u_j)\leq \alpha_n-x-L, \mathcal{H}_{\zeta_n, n}^u \bigg)\nonumber\\
	&\leq \P(\exists u\in \mathbb{T}: V(u)\leq -x-L/2 )\leq e^{-x-L/2}\lesssim (L+1)^{-1}e^{-x}. 
\end{align}
If $L< 2\alpha_n$, then for $n$ large enough $L\leq n^{(2b+1)/4}$.  Therefore, let $L_1$ be the unique integer such that $-L_1/2-1 < \alpha_n -L_1\leq -L_1/2$, then combining \eqref{Eq1} and \eqref{Step28} ,
\begin{align}
	& \P\bigg(\exists u\in \mathbb{T}_n, V(u)\leq \alpha_n -x, \min_{\tau_{\zeta_n}^{(u)}\leq j\leq n} V(u_j)\leq \alpha_n-x-L, \mathcal{H}_{\zeta_n, n}^u \bigg)\nonumber\\
	&\lesssim e^{-x-L_1/2}+ \sum_{j=[L]}^{L_1}  \P\bigg(\exists u\in \mathbb{T}_n, V(u)\leq \alpha_n -x, \min_{\tau_{\zeta_n}^{(u)}\leq j\leq n} V(u_j)-(\alpha_n -x)\in [-j, -j+1], \mathcal{H}_{\zeta_n, n}^u \bigg)\nonumber\\
	&\lesssim e^{-x-L_1/2}+ \sum_{j=[L]}^{L_1}  (1+j)^{-2} e^{-x}\lesssim (1+L)^{-1}e^{-x},
\end{align}
which implies \eqref{Eq2}. 

\hfill$\Box$

\begin{lemma}\label{lem:bad-event3}
	For any $\varepsilon, L>0$, there exists $T=T(\varepsilon, L)$ and $n_2(\varepsilon, L, T)>T$ such that for any $n\geq n_2$ and any $x\geq 0$, 
	\begin{align}
		& \P\bigg(\exists u\in\mathbb{T}_n, V(u)\leq \alpha_n-x, \min_{\tau_{\zeta_n}^{(u)}\leq j\leq n} V(u_j)\geq y_n, \tau_{\zeta_n}^{(u)}\leq n-T, \tau_{\zeta_n}^{(2,u)}>n\bigg)\nonumber\\
		&\leq \varepsilon e^{-x}.
	\end{align}
\end{lemma}
\textbf{Proof: } Combining the union bound and many-to-one formula \eqref{Many-to-one}, it holds that 
\begin{align}
	& \P\bigg(\exists u\in\mathbb{T}_n, V(u)\leq y_n+L, \min_{\tau_{\zeta_n}^{(u)}\leq j\leq n} V(u_j)\geq y_n, \tau_{\zeta_n}^{(u)}\leq n-T, \tau_{\zeta_n}^{(2,u)}>n\bigg)\nonumber\\
	&\leq  \E\bigg( \sum_{|u|=n} 1_{\left\{ V(u)\leq y_n+L, \min_{\tau_{\zeta_n}^{(u)}\leq j\leq n} V(u_j)\geq y_n, \tau_{\zeta_n}^{(u)}\leq n-T, \tau_{\zeta_n}^{(2,u)}>n \right\}} \bigg)\nonumber\\
	& = \mathbf{E}\left(e^{S_n}  1_{\left\{ S_n \leq y_n+L, \min_{\tau_{\zeta_n}\leq j\leq n} S_j \geq y_n, \tau_{\zeta_n} \leq n-T, \tau_{\zeta_n}^{(2)}>n \right\}} \right). 
\end{align}
By \eqref{Step32}, the above inequality has upper bound
\begin{align}\label{Step29}
	& \P\bigg(\exists u\in\mathbb{T}_n, V(u)\leq y_n+L, \min_{\tau_{\zeta_n}^{(u)}\leq j\leq n} V(u_j)\geq y_n, \tau_{\zeta_n}^{(u)}\leq n-T, \tau_{\zeta_n}^{(2,u)}>n\bigg)\nonumber\\
	& \leq  \sum_{k=1}^{[L]+1} e^{y_n+k}\mathbf{P}\left(S_n-y_n\in [k-1,k], \min_{\tau_{\zeta_n}\leq j\leq n} S_j \geq y_n, \tau_{\zeta_n}\leq n-T, \tau_{\zeta_n}^{(2)}>n \right)\nonumber\\
	&\lesssim e^{-x} \sum_{k=1}^{[L]+1}e^{k-L}  \left(\sum_{i=T}^{\infty} \mathbf{P}\left(\max_{0\leq j\leq i} S_j\leq L+1\right) +\frac{1}{n} \right)\nonumber\\
	& \lesssim e^{-x} \left( \sum_{i=T}^{\infty} \mathbf{P}\left(\max_{0\leq j\leq i} S_j\leq L+1\right) +\frac{1}{n} \right).
\end{align}
Recall the definition of the renewal function $R$ in \eqref{Def-of-renewal-function},  $\lim_{T\to\infty} \sum_{i=T}^\infty  \mathbf{P}\left(\max_{j\leq i}S_j \leq L+1\right)=0$, which together with \eqref{Step29} implies the desired result.

\hfill$\Box$

\begin{prop}\label{prop:Tightness}
	For any $x\in \R$,
	\begin{align}
		\P(M_n\leq \alpha_n -x)\lesssim e^{-x}.
	\end{align}
\end{prop}
\textbf{Proof: } 
Since $\P(M_n\leq \alpha_n -x)\leq 1\leq  e^{-x}$ when $x\leq 0$, we only consider the case $x\geq 0$.
From Lemma \ref{lem:bad-event1}, we see that 
\begin{align}\label{Step30}
	\P(M_n\leq \alpha_n -x)& \leq  \P(\exists u\in \mathbb{T}_n, V(u)\leq \alpha_n -x, \tau_{\zeta_n}^{(2,u)}\leq n ) \nonumber\\
	&\quad +  \P(\exists u\in \mathbb{T}_n, V(u)\leq \alpha_n -x, \tau_{\zeta_n}^{(u)}> n ) \nonumber\\
	&\quad + \P\left(\exists u\in \mathbb{T}_n, V(u)\leq \alpha_n -x, \mathcal{H}_{\zeta_n, n}^u\right) \nonumber\\
	&\lesssim e^{-x}+ \P\left(\exists u\in \mathbb{T}_n, V(u)\leq \alpha_n -x,\mathcal{H}_{\zeta_n, n}^u\right),
\end{align}
where recall the definition of $\mathcal{H}_{\zeta_n, n}^u$ in \eqref{Def-of-y}. 
Combining Lemma \ref{lem:bad-event2} and \eqref{Step30}, fix any $L\in \N$, when $n$ is large enough,
\begin{align}
	& \P(M_n\leq \alpha_n -x)\nonumber\\
	&\lesssim_L e^{-x}+  \P\bigg(\exists u\in \mathbb{T}_n, V(u)\leq \alpha_n -x, \min_{\tau_{\zeta_n}^{(u)}\leq j\leq n} V(u_j) >y_n,  \mathcal{H}_{\zeta_n, n}^u  \bigg)\nonumber\\
	&\leq e^{-x}+  \E\bigg(\sum_{|u|=n} 1_{\left\{ V(u)\leq y_n+L, \min_{\tau_{\zeta_n}^{(u)}\leq j\leq n} V(u_j) >y_n, \mathcal{H}_{\zeta_n, n}^u \right\}} \bigg).
\end{align}
Now it follows from the many-to-one formula \eqref{Many-to-one} and Lemma \ref{lem3} (ii) that 
\begin{align}
	\P(M_n\leq \alpha_n -x) & \lesssim_L  e^{-x}+  \mathbf{E}\left(e^{S_n}1_{\left\{ S_n\leq y_n+L, \min_{\tau_{\zeta_n}\leq j\leq n} S_j >y_n, \tau_{\zeta_n} \leq n<  \tau_{\zeta_n}^{(2)} \right\}}\right)\nonumber\\
	& \leq e^{-x}+  e^{\alpha_n -x}\mathbf{P}\left( S_n\leq y_n+L, \min_{\tau_{\zeta_n}\leq j\leq n} S_j >y_n, \tau_{\zeta_n} \leq n<  \tau_{\zeta_n}^{(2)} \right)\nonumber\\
	& \lesssim e^{-x}+  e^{\alpha_n -x} \sum_{i=1}^n \sum_{k=1}^{[L]+1}\left( \mathbf{P}\left(\max_{j\leq n-i} S_j \leq k\right)  +\frac{1}{n^2} \right)e^{-\alpha_n}\nonumber\\
	&\lesssim e^{-x} +(L+1)R(L+1)e^{-x},
\end{align}
as desired. 

\hfill$\Box$

Recall the definition of $\mathcal{S}$ in \eqref{Def-of-test-function}.
The following result is a refinement of Proposition \ref{prop:Tightness}, whose proof is postponed to Section \ref{S4}. 

\begin{prop}\label{prop:key-proposition}
	Suppose that $f\in \mathcal{S}$ for some $R_f>0 $. 
	Then for any $\varepsilon>0$, there exists a constant $A=A(\varepsilon, R_f)>0$ and $N=N(\varepsilon, R_f)\in \N $ such that for all $n\geq N$ and $x\in [A, n^{(3-2b)/4}]$, 
	\begin{align}\label{First-Ineq}
		\left|\E \left(1- e^{-\sum_{|u|=n} f(V(u)- (\alpha_n -x))} 1_{\left\{ M_n> \alpha_n -x \right\}}\right)- C^*(f)e^{-x}\right|\leq \varepsilon e^{-x},
	\end{align}
	where recall that $C^*(f)$ is given in \eqref{Def-of-C-star-f}. Moreover,
	\begin{align}\label{Second-Ineq}
		\left|\E \left(1- e^{-\sum_{|u|=n} f(V(u)- (\alpha_n -x))} \right)- \left(C^*(f)- C^*(0)\right)e^{-x}\right|\leq \varepsilon e^{-x}.
	\end{align}
\end{prop}

Note that Proposition \ref{prop:Tightness} and Proposition \ref{prop:key-proposition} imply the finiteness of $C^*(f)$. Now we are ready to prove Theorems \ref{thm1} and \ref{thm2}.

\textbf{Proof of Theorem \ref{thm1}: }  We just prove \eqref{Joint-prob} here, the proof of \eqref{Extremal-process} is similar with \eqref{First-Ineq} replaced by \eqref{Second-Ineq}. 
By Proposition \ref{prop:Tightness}, we see that 
\begin{align}
	\sum_{n=1}^\infty \P(M_n < \alpha_n - 2\log n) \lesssim \sum_{n=1}^\infty e^{-2\log n}<\infty,
\end{align}
which implies that almost surely,
\begin{align}\label{M-n-infinite}
	\liminf_{n\to\infty} \frac{M_n}{\lambda (mn)^b}\geq 1\quad \Longrightarrow\quad \lim_{n\to\infty} M_n=\infty. 
\end{align}
For each $K\geq 0$, define 
\begin{align}
	\mathcal{Z}(K):= \left\{u\in \mathbb{T}:  V(u)\geq K, V(u_k)< K,\forall\ k<|u| \right\}.
\end{align}
According to \cite[Theorem 9]{Ky2000}, almost surely,
\begin{align}
	\lim_{K\to\infty} \sum_{u\in \mathcal{Z}(K)} e^{-V(u)}=W_\infty. 
\end{align}
For each fixed $x\in \R$ and $0< \varepsilon<\min \{ C^*(0), 1\}$, let $A=A(\varepsilon, R_f)$ and $N=N(\varepsilon, R_f)$ be the constant given in Proposition \ref{prop:key-proposition}. Now we choose $K$ sufficiently large such that  $K>A+1+2|x|$. Since almost surely $\# \mathcal{Z}(K)$ is finite (under the assumption that $\nu$ is finite and \eqref{M-n-infinite}), there exists $N_1= N_1(\varepsilon, A, K)\in \N$ such that for any $n\geq N_1$, 
\begin{align}
	\frac{1}{2}n^{(3-2b)/4} +1 \leq  (n-n^b)^{(3-2b)/4},\quad	\max_{k\in [n-n^b, n]} |\alpha_n - \alpha_{k}| \leq \varepsilon, \quad \mbox{and }\quad \P(G_{K,n})\geq 1-\varepsilon,
\end{align}
where 
\begin{align}
	G_{K,n}:= \left\{V(u)-x\leq n^{(3-2b)/4}/2,\ \forall u\in \mathcal{Z}(K) \right\}\cap \left\{\max\left\{|u|: u\in \mathcal{Z}(K)\right\} \leq n^b \right\}. 
\end{align}
Also, according to our construction of $K$,  for any $u\in \mathcal{Z}(K)$ and any $k\in [n-n^b, n]$,
\begin{align}
	t:= V(u)-x+\alpha_k-\alpha_n\geq K - |x|-  \varepsilon > A
\end{align}
and that
\begin{align}
	t= V(u)-x+\alpha_k-\alpha_n\leq n^{(3-2b)/4}/2+ \varepsilon \leq k^{(3-2b)/4}. 
\end{align}
Therefore, combining \eqref{First-Ineq} in Proposition \ref{prop:key-proposition} and the strong Markov property, recall that $f_x= f(\cdot +x)$, we have the lower bound
\begin{align}\label{Step72}
	& \E \left(e^{-\sum_{|u|=n} f(V(u)- \alpha_n )} 1_{\left\{ M_n> \alpha_n + x \right\}}\right) = \E \left(e^{-\sum_{|u|=n} f_x(V(u)- \alpha_n -x)} 1_{\left\{ M_n> \alpha_n + x \right\}}\right)  \nonumber\\
	&  \geq \E \left(e^{-\sum_{|u|=n} f_x(V(u)- \alpha_n -x)} 1_{\left\{ M_n> \alpha_n + x \right\}}1_{G_{K,n}}\right)  \nonumber\\
	& =  \E\bigg(1_{G_{K,n}}\prod_{u\in \mathcal{Z}(K)}\E  \left(e^{-\sum_{|v|=k} f_x(V(v)- \alpha_k +t)} 1_{\left\{ M_k> \alpha_k -t \right\}}\right)\Big|_{k=n-|u|, t= V(u)-x+\alpha_k-\alpha_n} \bigg). 
\end{align} 
Applying Proposition \ref{prop:key-proposition} and noticing that $k\geq n-n^b\geq N$ and $t\in [A, k^{(3-2b)/4}]$, we deduce from \eqref{Step72} that
\begin{align}\label{Step72'}
	& \E \left(e^{-\sum_{|u|=n} f(V(u)- \alpha_n )} 1_{\left\{ M_n> \alpha_n + x \right\}}\right)  \nonumber\\
	&  \geq  \E\bigg(1_{G_{K,n}}\prod_{u\in \mathcal{Z}(K)}\left( 1-(C^*(f_x)+\varepsilon)e^{x-V(u)+\alpha_n -\alpha_k}\right) \bigg)\nonumber\\
	&\geq \E\bigg(\prod_{u\in \mathcal{Z}(K)}\left( 1-(C^*(f_x)+\varepsilon)e^{x-V(u)+\varepsilon}\right) \bigg)-\varepsilon.
\end{align}
Since $\max\left\{e^{-V(u)}: u\in\mathcal{Z}(K)\right\}\leq e^{-K}\to 0$ as $K\to\infty$, using the fact that $1-z\sim e^{-z}$ as $z\to 0$, we have 
$\prod_{u\in \mathcal{Z}(K)} \left( 1-(C^*(f_x)+\varepsilon) e^{x-V(u)+ \varepsilon}  \right)\to \exp\big\{ -(C^*(f_x)+\varepsilon) e^{x+\varepsilon}W_\infty \big\} $, which together with \eqref{Step72'} implies that 
\begin{align}
	\liminf_{n\to\infty}  \E \left(e^{-\sum_{|u|=n} f(V(u)- \alpha_n )} 1_{\left\{ M_n> \alpha_n + x \right\}}\right)   & \geq \liminf_{K\to\infty} \E\bigg(\prod_{u\in \mathcal{Z}(K)} \left(1-(C^*(f_x)+\varepsilon) e^{x-V(u)+ \varepsilon} \right)  \bigg)	 -\varepsilon\nonumber\\
	& = \E\left( \exp\left\{-(C^*(f_x)+\varepsilon) e^{x+\varepsilon}W_\infty \right\}\right)-\varepsilon.
\end{align}
Taking $\varepsilon\to 0$ in the above inequality yields that 
\begin{align}
	\liminf_{n\to\infty}  \E \left(e^{-\sum_{|u|=n} f(V(u)- \alpha_n )} 1_{\left\{ M_n> \alpha_n + x \right\}}\right) & \geq  \E\left( \exp\left\{- C^*(f_x) e^{x}W_\infty \right\}\right). 
\end{align}
The proof for the upper bound is similar and we omit the details here, this completes the proof of the theorem. 

\hfill$\Box$

\textbf{Proof of Theorem \ref{thm2}: } We firstly show that for any non-negative continuous function $f$ with bounded support, 
\begin{align}\label{Goal1}
	\E \left(e^{-\int f(y) \mathcal{E}_\infty(\mathrm{d}y) }\right) = \E\left( \exp\left\{- \left(C^*(f)-C^*(0)\right)W_\infty \right\}\right). 
\end{align}
Since $\mathbb{T}^{(i)}$ is independent of $\{(p_i, q_i), i\in \mathbb{N}\}$ and $W_\infty$, it holds that
\begin{align*}
	\E \left(e^{-\int f(y) \mathcal{E}_\infty(\mathrm{d}y) } \big| (p_i, q_i), i\in\mathbb{N}, W_\infty\right)= \prod_ {i=1}^\infty \mathbb{E}\bigg(\exp\bigg\{-1_{\{M_j\geq z\}}\sum_{|u|=j} f(V(u)-z)\bigg\}\bigg)\bigg|_{(z,j)=(p_i, q_i)}.
\end{align*}
Taking expectation with respect to $(p_i, q_i)$ in the above equation, we get that 
\begin{align}
	&\E \left(e^{-\int f(y) \mathcal{E}_\infty(\mathrm{d}y) } \big|  W_\infty\right)= \exp\bigg\{-W_\infty \ell_\infty m^a\sum_{j=0}^\infty \int_{\mathbb{R}} e^{-z}\mathbb{E}\left(1-  e^{-1_{\{M_j\geq z\}}\sum_{|u|=j} f(V(u)-z)}\right)\mathrm{d}z\bigg\}\nonumber\\
	& = \exp\bigg\{-W_\infty \ell_\infty m^a\sum_{j=0}^\infty \int_{\mathbb{R}} e^{-z}\mathbb{E}\left(\left(1-  e^{-\sum_{|u|=j} f(V(u)-z)}\right)1_{\{M_j\geq z\}}\right)\mathrm{d}z\bigg\}\nonumber\\
	&\stackrel{z:=  M_j-\zeta}{=}\exp\bigg\{-W_\infty \ell_\infty m^a\mathbb{E}\bigg(e^{-M_j}\sum_{j=0}^\infty \int_{0}^\infty  e^{\zeta}\left(1-  e^{-\sum_{|u|=j} f(V(u)-M_j+\zeta)}\right)\mathrm{d}\zeta\bigg)\bigg\}\nonumber\\
	& = \exp\left\{-\left(C^*(f)-C^*(0)W_\infty \right)\right\}. 
\end{align}
Now taking expectation with respect to $W_\infty$ implies \eqref{Goal1}. Therefore, by Theorem \ref{thm1},  it suffices to prove that $\mathcal{E}_\infty$ is almost surely a local finite measure, which is equivalent to show that $\mathbb{P}(\int f(y) \mathcal{E}_\infty(\mathrm{d}y)<\infty )=1$ for any non-negative continuous function $f$ with bounded support. From \eqref{Goal1}, we see that 
\begin{align}\label{Step79}
	&\mathbb{P}\left(\int f(y) \mathcal{E}_\infty(\mathrm{d}y)<\infty\right) \nonumber\\
	&= \lim_{\theta\downarrow 0} \E \left(e^{-\theta \int f(y) \mathcal{E}_\infty(\mathrm{d}y) }\right) = \lim_{\theta \downarrow 0}\E\left( \exp\left\{- \left(C^*(\theta f)-C^*(0)\right)W_\infty \right\}\right). 
\end{align}
According to dominated convergence theorem, it is easy to verify from \eqref{Def-of-C-star-f} that 
\begin{align}\label{Step80}
	&\lim_{\theta\downarrow 0} \left(C^*(\theta f)-C^*(0)\right) \nonumber\\
	&= \ell_\infty m^a   \lim_{\theta\downarrow 0}   \sum_{j=0}^\infty \E  \left(  e^{-M_j}   \int_0^\infty e^{z} \left(1- e^{-\sum_{|v|=j} \theta f(V(v)-M_j+z)} \right) \mathrm{d} z\right) =0.
\end{align}
Combining \eqref{Step79} and \eqref{Step80}, we complete the proof of the theorem. 

\hfill$\Box$

\section{Proof of Theorem \ref{thm3}}\label{S5}

We prove Theorem \ref{thm3} in this section, inspired by the arguments in \cite{AHS}. We first gather some useful facts which will be used in the proof. 

Combining \eqref{Step1} and Lemma \ref{lem1} (ii) with $x:= \frac{1}{2}mn$, we get
\begin{align}
	\mathbf{P}\left(S_n\leq \frac{1}{2}mn\right) &\leq 	\mathbf{P}\left(S_n\leq \frac{1}{2}mn, \min_{1\leq j\leq n} X_j\geq -\zeta_n\right)+ n\mathbf{P}(X\leq -\zeta_n)\nonumber\\
	& \lesssim e^{-\theta_n mn/2} + n \zeta_n^{a+1-b} e^{-\lambda \zeta_n^b}.
\end{align}
Since $\zeta_n\asymp n$ and that $\theta_n mn/2 \asymp n^{b}$, there exists a constant $C_{up}>0$ such that for large $n$,
\begin{align}\label{upp1}
	\mathbf{P}\left(S_n\leq \frac{1}{2}mn\right) \leq e^{-C_{up}n^b}.
\end{align}
For the lower bound, since $p_0<1$, by \eqref{Equal-Conditions}, for any $-x\leq x_0$, 
\begin{align}
	& \P\left(\exists u \in \mathbb{T}_1: V(u)\leq -x\right)\nonumber\\
	& \geq \P\left(Y_1\in (-x-1, -x), \nu>0\right)= (1-p_0)\int_{-x-1}^x \ell(y)|y|^a e^{-\lambda |y|^b +y}\mathrm{d}y.
\end{align}
Therefore, there exists a constant $C_{low}>0$ such that for large $x>0$,
\begin{align}\label{low1}
	\P\left(\exists u\in \mathbb{T}_1: V(u)\leq -x\right)\geq e^{-x}e^{-C_{low}x^b}. 
\end{align}

Define $\P^*(\cdot):= \P(\cdot|  \# \mathbb{T}=\infty)$,  $a_*:= 	\inf_{\Theta\in \partial\mathbb{T}}\limsup_{n\to\infty} \frac{V(\Theta_n)}{n^{\frac{1}{2-b}}} $ and set $a_*:=0$ if $\partial \mathbb{T}=\emptyset$.  Denote $N_n:= \# \mathbb{T}_n$ by the number of the particles alive at generation $n$. Recall that $\mathcal{F}_n$ is the natually filtration of the branching random walk, then by the branching property, for any $a>0$,
\begin{align}
	\mathbb{P}\left(a_* \geq a\Big| \mathcal{F}_n\right) = \P(a_*\geq a)^{N_n}.
\end{align}
Therefore, $\P(a_*\geq a)^{N_n}$ is a non-negative martingale and thus converges $\P$-almost surely to $1_{\{a_*\geq a\}}$. Noticing that this process survives with positive probability and $N_n\to \infty$ $\P^*$-a.s., we obtain that $\P^*$-a.s., $1_{\{a_*\geq a\}}= \lim_{n\to\infty} \P(a_*\geq a)^{N_n}= 1_{\{ \P(a_*\geq a) =1\}},$ which implies that $a_*$ is a constant. Therefore, it remains to show that $a_*\in (0,\infty)$. 

\paragraph{Proof of Theorem \ref{thm3} for $a_*>0$.}  
Fix a sufficiently small $a>0$ such that $2a^2 < C_{up}m (1-b)$. It suffices to show that for any $n_0\geq 1$ and $ K>0$,
\begin{align}\label{Goal3}
	\lim_{n\to\infty} \mathbb{P}\left(\exists u\in \mathbb{T}_n:\ -K\leq V(u_\ell)\leq a \ell^{\frac{1}{2-b}},\ \forall \ n_0\leq \ell \leq n \right)=0. 
\end{align}
Indeed, by \eqref{M-n-infinite}, we can choose a suitable $K>0$ such that $\P^*(\inf_{u\in \mathbb{T}} V(u)>-K)>0$, which also implies that $\P(\inf_{u\in \mathbb{T}} V(u)>-K)>0$ since $\P^*(\inf_{u\in \mathbb{T}} V(u)>-K)\leq \frac{1}{\P(\# \mathbb{T}=\infty)} \P(\inf_{u\in \mathbb{T}} V(u)>-K)$. Then from \eqref{Goal3},  on the event $\left\{ \inf_{u\in \mathbb{T}} V(u)>-K\right\}$, $\P$-almost surely for all $n_0$, there is no $\Theta\in \partial \mathbb{T}$ such that $V(\Theta_\ell)\leq a \ell^{\frac{1}{2-b}} $ for all $\ell \geq n_0$. Therefore, $a_*\geq a$ holds with positive $\P^*$-probability, which also implies that $a_* \geq a>0$. 

Now in the rest part of the proof, we aim to prove \eqref{Goal3}.  
According to our choice of $a$,  $\frac{2a}{m}<\frac{C_{up}(1-b)}{a}$. Therefore,  there exists some constant $C$ such that
\begin{align}\label{Choice-of-C}
	\frac{2a}{m}< C^{\frac{1-b}{2-b}}<\frac{C_{up}(1-b)}{a}\ \Longleftrightarrow\ \frac{1}{2}m C> a C^{\frac{1}{2-b}}\ \mbox{and}\ a C^{\frac{1}{2-b}}< C_{up}C^b(1-b).
\end{align}
Define $r_i:= n_0+[C i^{\frac{2-b}{1-b}}]$. For each $L\in \mathbb{N}$, by Markov's inequality and many-to-one formula \eqref{Many-to-one},
\begin{align}
	& \mathbb{P}\left(\exists u\in \mathbb{T}_{r_L}:\ -K\leq V(u_\ell)\leq a \ell^{\frac{1}{2-b}},\ \forall \ n_0\leq \ell \leq r_L \right)\nonumber\\
	&\leq \mathbb{E}\bigg(\sum_{|u|\in r_L} 1_{\left\{  -K\leq V(u_{r_i})\leq a r_i^{\frac{1}{2-b}},\ \forall \ 1\leq i \leq L \right\}} \bigg) = \mathbf{E}\bigg(e^{S_{r_L}}1_{\left\{ -K\leq S_{r_i}\leq a r_i^{\frac{1}{2-b}},\ \forall \ 1\leq i \leq L \right\}} \bigg)\nonumber\\
	&\leq e^{a r_L^{\frac{1}{2-b}}} \mathbf{P}\left( -K\leq S_{r_i}\leq a r_i^{\frac{1}{2-b}},\ \forall \ 1\leq i \leq L \right) \leq e^{a r_L^{\frac{1}{2-b}}} \prod_{i=1}^{L}\mathbf{P}\left( S_{r_i-r_{i-1}}\leq a r_i^{\frac{1}{2-b}}+K\right) .
\end{align}
We use $a_n\sim b_n$ to denote $\lim_{n\to\infty} a_n/b_n=1$. 
Combining \eqref{Choice-of-C},  $ar_i^{\frac{1}{2-b}}+K \sim a C^{\frac{1}{2-b}} i^{\frac{1}{1-b}}$ and $r_{i}-r_{i-1}\sim Ci^{\frac{1}{1-b}}$, we can choose $N_0$ large enough such that $ar_i^{\frac{1}{2-b}}+K < \frac{1}{2}m (r_{i}-r_{i-1})$ and that \eqref{upp1} holds for all $i\geq N_0$. Now for $L\geq N_0$ large enough, the above probability has upper bound
\begin{align}
	& \mathbb{P}\left(\exists u\in \mathbb{T}_{r_L}:\ -K\leq V(u_\ell)\leq a \ell^{\frac{1}{2-b}},\ \forall \ n_0\leq \ell \leq r_L \right)\nonumber\\
	& \leq e^{a r_L^{\frac{1}{2-b}}} \prod_{i=N_0}^{L}\mathbf{P}\left( S_{r_i-r_{i-1}}\leq  \frac{1}{2}m (r_{i}-r_{i-1})\right) \leq \exp\bigg\{a r_L^{\frac{1}{2-b}} - C_{up}\sum_{i=N_0}^L (r_{i}-r_{i-1})^b\bigg\}.
\end{align}
Noticing that $\sum_{i=N_0}^L (r_{i}-r_{i-1})^b\sim C^b \sum_{i=N_0}^L i^{\frac{b}{1-b}}\sim C^b(1-b) L^{\frac{1}{1-b}}$ and that $ar_L^{\frac{1}{2-b}}\sim aC^{\frac{1}{2-b}} L^{\frac{1}{1-b}}$, it follows from \eqref{Choice-of-C} that 
\begin{align}
	\lim_{L\to\infty} \mathbb{P}\left(\exists u\in \mathbb{T}_{r_L}:\ -K\leq V(u_\ell)\leq a \ell^{\frac{1}{2-b}},\ \forall \ n_0\leq \ell \leq r_L \right) =0,
\end{align}
which implies \eqref{Goal3}. We are done.

\hfill$\Box$

\paragraph{Proof of Theorem \ref{thm3} for $a_*<\infty$.}  Let $r_1$ be a fixed large integer and define $r_i:= r_1+ [(i-1)^{\frac{2-b}{1-b}}]$. We use $u<v$ to denote that $v$ is an ancestor of $u$.
Let $\mathcal{Z}_1:= \mathbb{T}_{r_1}$ be the set of particles alive at time $r_1$ and for each $C>0$ and $i\geq 2$, define $\mathcal{Z}_i$ recursively by
\begin{align}
	\mathcal{Z}_i:= &\cup_{v\in \mathcal{Z}_{i-1}} \Big\{ u\in \mathbb{T}_{r_i}, u<v:\ V(u_k)-V(v)\leq C(k-r_{i-1}),\forall  r_{i-1}\leq k\leq r_i,\nonumber\\
	& \quad V(u_{r_i-1})-V(v)\geq 2C i^{\frac{b}{1-b}}, \ V(u)-V(v)\leq Ci^{\frac{b}{1-b}}
	\Big\}.
\end{align}
It is easy to see that $\mathcal{Z}_i\subset \mathbb{T}$ and that $(\# \mathcal{Z}_i)$ is an time-inhomogeneous Galton-Watson tree. If $\mathcal{Z}_i$ survives with positive probability, then there exists an infinite ray $\Theta$ of $( \mathcal{Z}_i)\subset \mathbb{T}$ such that 
$V(\Theta_{r_i})-V(\Theta_{r_{i-1}})\leq Ci^{\frac{b}{1-b}}$ for all $i\in \mathbb{N}$, which implies that 
$V(\Theta_{r_i})\leq V(\Theta_{r_1})+C \sum_{j=2}^i j^{\frac{b}{1-b}}\leq V(\Theta_{r_1})+ Ci^{\frac{1}{1-b}}$. Moreover, for any $r_i\leq n< r_{i+1}$, we also have that $V(\Theta_n)\leq V(r_{i})+C(n-r_i)\leq V(\Theta_{r_1})+ Ci^{\frac{1}{1-b}}+ C(r_{i+1}-r_i)$. Since $r_{i+1}-r_i\sim i^{\frac{1}{1-b}}\sim r_i^{\frac{1}{2-b}}$, we conclude that when $n$ is large enough, $V(\Theta_n) \leq 3Cn^{\frac{1}{2-b}}$, which implies that $a_* \leq 3C$.

To show that $(\# \mathcal{Z}_i)$ survives with positive probability, define $b_i:= r_i-r_{i-1}-1$ and 
\[
\nu_i:= \sum_{|w|=b_i} 1_{\left\{V(w_k)\leq Ck, \forall \ k\leq b_i, V(w)\geq 2C i^{\frac{b}{1-b}}\right\}} 1_{\left\{\exists u\in \mathbb{T}_{b_i+1}, u<w: V(u)\leq Ci^{\frac{b}{1-b}} \right\}}.
\]
Following the same argument as \cite{AHS}, to prove $a_*<\infty$, it remains to show that for any $K_0\geq 1$,
\begin{align}\label{Goal4}
	\limsup_{C\to\infty} \limsup_{i\to\infty} \P(\nu_i \leq K_0)\leq 1-\mathbb{P}(\# \mathbb{T}=\infty).
\end{align}
It follows from the Markov property that 
\begin{align}
	\E\left(\nu_i \big|\mathcal{F}_{b_i}\right)= \sum_{|w|=b_i} 1_{\left\{V(w_k)\leq Ck, \forall \ k\leq b_i, V(w)\geq 2C i^{\frac{b}{1-b}}\right\}} \P\left( \exists u\in \mathbb{T}_1: V(u)\leq -x\right)\big|_{x= V(w)-Ci^{\frac{b}{1-b}}}.
\end{align}
Since $Ci^{\frac{b}{1-b}} \leq V(w)- Ci^{\frac{b}{1-b}} \leq Cb_i$, by \eqref{low1}, when $i$ is large enough, we have that 
\begin{align}\label{Step81}
	\E\left(\nu_i \big|\mathcal{F}_{b_i}\right) &\geq e^{Ci^{\frac{b}{1-b}}- C_{low}(Cb_i)^b}  \sum_{|w|=b_i} 1_{\left\{V(w_k)\leq Ck, \forall \ k\leq b_i, V(w)\geq 2C i^{\frac{b}{1-b}}\right\}} e^{-V(w)}\nonumber\\
	&=: e^{Ci^{\frac{b}{1-b}}- C_{low}(Cb_i)^b}  W_i'.
\end{align}
Define $p(w):=  \P\left( \exists u\in \mathbb{T}_1: V(u)\leq -x\right)\big|_{x= V(w)-Ci^{\frac{b}{1-b}}}.$ Then combining the inequality
\begin{align}
	\mbox{Var}\left(\nu_i \big| \mathcal{F}_{b_i}\right)=  \sum_{|w|=b_i} 1_{\left\{V(w_k)\leq Ck, \forall \ k\leq b_i, V(w)\geq 2C i^{\frac{b}{1-b}}\right\}} p(w)(1-p(w))\leq \E\left(\nu_i \big|\mathcal{F}_{b_i}\right),
\end{align}
\eqref{Step81} and Markov's inequality,  for any $\varepsilon>0$,
\begin{align}\label{Step82}
	&\limsup_{i\to\infty} \P\left(\nu_i- \E\left(\nu_i \big|\mathcal{F}_{b_i}\right) \leq -\varepsilon\E\left(\nu_i \big|\mathcal{F}_{b_i}\right), W_i'\geq \varepsilon \right)\nonumber\\
	&\leq \limsup_{i\to\infty} \E\left[ \frac{ \mbox{Var}\left(\nu_i \big| \mathcal{F}_{b_i}\right)}{\varepsilon^2 \left( \E\left(\nu_i \big|\mathcal{F}_{b_i}\right) \right)^2}1_{\{W_i'\geq \varepsilon\}}\right] \leq \limsup_{i\to\infty} \E\left[ \frac{ 1}{\varepsilon^2  \E\left(\nu_i \big|\mathcal{F}_{b_i}\right) }1_{\{W_i'\geq \varepsilon\}}\right] \nonumber\\
	& \leq \frac{1}{\varepsilon^3}\limsup_{i\to\infty}e^{C_{low}(Cb_i)^b-Ci^{\frac{b}{1-b}} }.
\end{align}
Since $(Cb_i)^b\sim C^bi^{\frac{b}{1-b}}$, let $C$ be large enough such that $C>C_{low} C^b$, then from \eqref{Step81}, on the event $W_i'\geq \varepsilon$, for large $i$, it holds that $(1-\varepsilon) \E\left(\nu_i \big|\mathcal{F}_{b_i}\right)  \geq (1-\varepsilon)\varepsilon e^{Ci^{\frac{b}{1-b}}- C_{low}(Cb_i)^b}\geq K_0$. Therefore, we conclude from \eqref{Step82} that
\begin{align}\label{Step83}
	\limsup_{i\to\infty} \P(\nu_i \leq K_0)& \leq  \limsup_{i\to\infty} \P(W_i'\leq \varepsilon) + \limsup_{i\to\infty} \P\left(\nu_i- \E\left(\nu_i \big|\mathcal{F}_{b_i}\right) \leq -\varepsilon\E\left(\nu_i \big|\mathcal{F}_{b_i}\right), W_i'\geq \varepsilon \right)\nonumber\\
	&=  \limsup_{i\to\infty} \P(W_i'\leq \varepsilon).
\end{align}
Recall that $W_n = \sum_{|u|=n}e^{-V(u)}$ is a non-negative martingale with $L^1(\P)$ limit $W_\infty$.  According to the many-to-one formula \eqref{Many-to-one},
\begin{align}\label{Step84}
	\E\left(W_i -W_i'\right)&= 1- \mathbf{P}\left(S_k\leq Ck, \forall \ k\leq b_i,\ S_{b_i}\geq 2 Ci^{\frac{b}{1-b}}\right)\nonumber\\
	&\leq \mathbf{P}\left(\max_{k\geq 1} \frac{S_k}{k}>C\right)+ \mathbf{P}(S_{b_i}> 2 Ci^{\frac{b}{1-b}}).
\end{align}
Since $b_i\sim i^{\frac{1}{1-b}}$, the right-hand side of the above inequality tends to $0$ as $C\to\infty$. Therefore, combining Markov's inequality,  \eqref{Step83} and \eqref{Step84},
\begin{align}\label{Step85}
	\limsup_{C\to\infty} \limsup_{i\to\infty} \P(\nu_i \leq K_0)& \leq \limsup_{C\to\infty}\limsup_{i\to\infty} \left(\P(W_i\leq 2\varepsilon)+ \P(W_i-W_i'>\varepsilon)\right)\nonumber\\
	&\leq \limsup_{i\to\infty}  \P(W_i\leq 2\varepsilon) + \frac{1}{\varepsilon}\limsup_{C\to\infty}\limsup_{i\to\infty}\E\left(W_i -W_i'\right)\nonumber\\
	&=  \limsup_{i\to\infty}  \P(W_i\leq 2\varepsilon).
\end{align}
Noticing that $W_i$ converges almost surely to $W_\infty$ and that under \eqref{Assumption2},  $\{W_\infty>0\}$ coincides with $\{\# \mathbb{T}=\infty\}$, taking $\varepsilon\downarrow 0$ in \eqref{Step85}, we get \eqref{Goal4}. This implies the desired result. 

\hfill$\Box$

\section{Proof of Proposition \ref{prop:key-proposition}} \label{S4} 

We give the proof of Proposition \ref{prop:key-proposition} in this section.

For each $u\in \mathbb{T}\setminus \{\emptyset\}$, we define $\mathbb{B}(u)$ by the set of siblings of $u$. Recall that the probability measure $\Q$ is defined in \eqref{Change-of-measure} and that $\{w_n\}_{n\geq 0}$ is the spine process.  Let $B>0$ be a large constant and $J$ be a large integer which will be determined later. We say that a particle $u\in \mathbb{T}_n$ is a good vertex if for any $x\geq 0$, 
\begin{align}\label{Def-of-good-vertex}
	\tau_{\zeta_n}^{(2,u)}>n \geq \tau_{\zeta_n}^{(u)}>J\quad \mbox{and} \quad \sum_{v\in\mathbb{B}(u_k)} e^{-(V(v)+x)} \leq \begin{cases}
		e^{B-x},\quad & \mbox{if} \ 1\leq k\leq J;\\
		e^{-mk/3},\quad &\mbox{if} \ J<k<\tau_{\zeta_n}^{(u)}.
	\end{cases}
\end{align}
Recall the definition of $y_n$ in \eqref{Def-of-y}.

\begin{lemma}\label{lem:bad-even4}
	For any $\varepsilon, L>0, T\in\N$, there exists $J_0=J_0(\varepsilon, L, T)$ such that for any $J\geq J_0$, there exists $B_0=B_0(\varepsilon, L, T, J)$ such that for any $B\geq B_0$, when $n$ is large enough, for any $0\leq x\leq  n^{(3-2b)/4}$, 
	\begin{align}\label{Goal}
		& \Q\bigg(V(w_n)\leq \alpha_n-x, \min_{\tau_{\zeta_n}^{(w_n)}\leq j\leq n} V(w_j)\geq y_n, \tau_{\zeta_n}^{(w_n)}\in [n-T,n], w_n\ \mbox{not good}\bigg)\nonumber\\
		&\leq \varepsilon e^{-\alpha_n}. 
	\end{align}
\end{lemma}
\textbf{Proof: } Define 
\begin{align}\label{Def-of-F-n}
	F_n:=\bigg\{ V(w_n)\leq \alpha_n-x, \min_{\tau_{\zeta_n}^{(w_n)}\leq j\leq n} V(w_j)\geq y_n, \tau_{\zeta_n}^{(w_n)}\in [n-T,n]\bigg\}. 
\end{align}
According to Lemma \ref{lem3} (i), when $n$ is large enough,
\begin{align}\label{Proof-part1}
	\Q\left(F_n, \tau_{\zeta_n}^{(2,w_n)}\leq n\right)\leq\Q\left( \tau_{\zeta_n}^{(2,w_n)}\leq n\right)= \mathbf{P}\left(\tau_{\zeta_n}^{(2)}\leq n\right)\leq \frac{\varepsilon}{4}e^{-\alpha_n}.
\end{align}

\paragraph{Step 1.} In this step, we show that there exists a constant $\Gamma= \Gamma(\varepsilon, L, T)>0$ such that for large $n$,
\begin{align}\label{Proof-part2}
	\Q\bigg(F_n, \min_{1\leq j<\tau_{\zeta_n}^{(w_n)}} V(w_j)\leq -\Gamma,  \tau_{\zeta_n}^{(2,w_n)}> n\bigg)\leq \frac{\varepsilon}{4}e^{-\alpha_n}.
\end{align}
Noticing that by the Markov property, the left hand side of \eqref{Proof-part2} is equal to 
\begin{align}\label{Step37}
	& \Q\bigg(F_n, \min_{1\leq j<\tau_{\zeta_n}^{(w_n)}} V(w_j)\leq -\Gamma,  \tau_{\zeta_n}^{(2,w_n)}> n\bigg)\nonumber\\
	& = \sum_{i=n-T}^n \sum_{k=1}^{i-1} \mathbf{E}\bigg(1_{\left\{\min_{1\leq j\leq k-1}S_j>-\Gamma, S_k\leq -\Gamma, \min_{1\leq j\leq k} X_j\geq -\zeta_n   \right\}} \nonumber\\
	&\quad\quad\quad \times \mathbf{P}\left(S_{n-k}\leq \alpha_n-x-q, \min_{i-k\leq \ell\leq n-k} S_\ell \geq y_n-q, \tau_{\zeta_n}=i-k,\tau_{\zeta_n}^{(2)}> n-k\right)\bigg|_{q=S_k} \bigg)\nonumber\\
	&=:  \sum_{i=n-T}^n \sum_{k=1}^{i-1} \Upsilon_n(i,k).
\end{align}
When $k\geq n^{(3-2b)/4}$, then combining \eqref{Step16} and Markov's inequality, 
\begin{align}
	\Upsilon_n(i,k)& \leq \mathbf{P}(X_{i-k}<-\zeta_n) \mathbf{P}(S_k \leq -\Gamma) \leq \mathbf{P}(X <-\zeta_n) \mathbf{P}(|S_k- m k| \geq m k) \nonumber\\
	&\lesssim  e^{-\alpha_n} n^{1-b+ \lambda b m^{b-1} A_1} \frac{\mathbf{E}\left(|S_k-mk|^D\right)}{(mk)^D}.
\end{align}
From \eqref{Step17}, fixing any integer $D$ such that $D(3-2b)/8\geq 1-b+\lambda b m^{b-1}A_1 + 2$, we get 
\begin{align}\label{Step38}
	\sum_{i=n-T}^n \sum_{k=[n^{(3-2b)/4}]+1}^{i-1} \Upsilon_n(i,k)& \lesssim \sum_{i=n-T}^n \sum_{k=[n^{(3-2b)/4}]+1}^{i-1} e^{-\alpha_n} n^{1-b+ \lambda b m^{b-1} A_1} k^{-D/2}\nonumber\\
	&\leq (T+1) n \times e^{-\alpha_n} n^{1-b+ \lambda b m^{b-1} A_1} n^{-D(3-2b)/8}\nonumber\\
	&\leq \frac{T+1}{n}e^{-\alpha_n}. 
\end{align}
When $k\leq [n^{(3-2b)/4}]$, if $S_k \geq -n^{(3-2b)/4}$, then $y_n -S_k \leq 2n^{(3-2b)/4}\leq 4(n-k)^{(3-2b)/4}$ for large $n$. Therefore, combining Lemma \ref{lem3}(ii) and the fact that $\alpha_{n-k}= \alpha_n +o(1)$ uniformly in $1\leq k\leq [n^{(3-2b)/4}]$, we get
\begin{align}
	&\mathbf{P}\left(S_{n-k}\leq \alpha_n-x-q, \min_{i-k\leq \ell\leq n-k} S_\ell \geq y_n-q, \tau_{\zeta_n}=i-k,\tau_{\zeta_n}^{(2)}> n-k\right)\bigg|_{q=S_k}\nonumber\\
	& \lesssim \sum_{z=0}^{[L]} \left(\mathbf{P}\left(\max_{j\leq n-i+1} S_j \leq z+1\right) e^{-\alpha_{n-k}} +\frac{1}{n^2}e^{-\alpha_n}\right)\nonumber\\
	&\lesssim (L+1)e^{-\alpha_n},
\end{align}
which implies that 
\begin{align}\label{Step36}
	&\sum_{i=n-T}^n \sum_{k=1}^{[n^{(3-2b)/4}]} \Upsilon_n(i,k) \lesssim \sum_{i=n-T}^n \sum_{k=1}^{[n^{(3-2b)/4}]} \mathbf{P}(S_k< -n^{(3-2b)/4})\mathbf{P}(X_{i-k}<-\zeta_n)\nonumber\\
	&\quad + \sum_{i=n-T}^n \sum_{k=1}^{[n^{(3-2b)/4}]} \mathbf{P} \left(\min_{j\leq k-1}S_j>-\Gamma, S_k\leq -\Gamma \right) \times (L+1)e^{-\alpha_n}\nonumber\\
	& \leq (T+1)\mathbf{P}(X<-\zeta_n) \sum_{k=1}^{[n^{(3-2b)/4}]} \mathbf{P}(S_k< -n^{(3-2b)/4})\nonumber\\
	&\quad + (L+1)(T+1)e^{-\alpha_n} \mathbf{P}\left(\min_{k\geq 0} S_k \leq -\Gamma \right). 
\end{align}
According to Markov's inequality, we see that 
\begin{align}\label{Step35}
	\sum_{k=1}^{[n^{(3-2b)/4}]} \mathbf{P}(S_k< -n^{(3-2b)/4}) & \leq \sum_{k=1}^{[n^{(3-2b)/4}]} \mathbf{P}(S_k-mk< -n^{(3-2b)/4})\nonumber\\
	& \leq \sum_{k=1}^{[n^{(3-2b)/4}]} \frac{\mathbf{E}(|S_k-mk|^D)}{n^{(3-2b)D/4}}.
\end{align}
Combining \eqref{Step16}, \eqref{Step36} and \eqref{Step35}, fixing any integer $D$ such that $(3-2b)(D-2)/8 > 1-b+\lambda b m^{b-1}A_1 +1$, we obtain
\begin{align}\label{Step39}
	&\sum_{i=n-T}^n \sum_{k=1}^{[n^{(3-2b)/4}]} \Upsilon_n(i,k) \nonumber\\
	&\lesssim (T+1) e^{-\alpha_n} n^{1-b+ \lambda b m^{b-1} A_1} \sum_{k=1}^{[n^{(3-2b)/4}]} \frac{1}{n^{(3-2b)D/8}} + (L+1)(T+1)e^{-\alpha_n} \mathbf{P}\left(\min_{k\geq 0} S_k \leq -\Gamma \right)\nonumber\\
	&\lesssim \frac{T+1}{n}e^{-\alpha_n}+ (L+1)(T+1)e^{-\alpha_n} \mathbf{P}\left(\min_{k\geq 0} S_k \leq -\Gamma \right).
\end{align}
Therefore, combining \eqref{Step37}, \eqref{Step38} and \eqref{Step39}, we deduce that
\begin{align}
	& \Q\bigg(F_n, \min_{1\leq j<\tau_{\zeta_n}^{(w_n)}} V(w_j)\leq -\Gamma,  \tau_{\zeta_n}^{(2,w_n)}> n\bigg)\nonumber\\
	& \lesssim \frac{T+1}{n}e^{-\alpha_n}+ (L+1)(T+1)e^{-\alpha_n} \mathbf{P}\left(\min_{k\geq 0} S_k \leq -\Gamma \right),
\end{align}
which implies \eqref{Proof-part2} when both of $n$ and $\Gamma$ are  large enough. 

\paragraph{Step 2.} In this step, 
we prove that for $n$ large enough, 
\begin{align}\label{Proof-part3}
	\Q\bigg(F_n, \min_{1\leq j<\tau_{\zeta_n}^{(w_n)}} V(w_j)> -\Gamma,  \tau_{\zeta_n}^{(2,w_n)}> n, w_n\ \mbox{not good}\bigg)\leq \frac{\varepsilon}{2}e^{-\alpha_n}.
\end{align}
Once \eqref{Proof-part3} is proved, then we complete the proof of the lemma combining \eqref{Proof-part1}, \eqref{Proof-part2} and \eqref{Proof-part3}. 
For each $k$, define 
\[
\xi(w_k):=\sum_{v\in\mathbb{B}(w_k)} e^{-(V(v)-V(w_{k-1}))},
\]
then for each $k$, $\xi(w_k)$ are iid and that 
\[
\sum_{v\in \mathbb{B}(w_k)} e^{-(V(v)+x)}= e^{-(V(w_{k-1})+x)}\xi(w_k).
\]
Moreover, by the union bound, we have
\begin{align}
	& \Q\bigg(F_n, \min_{1\leq j<\tau_{\zeta_n}^{(w_n)}} V(w_j)> -\Gamma,  \tau_{\zeta_n}^{(2,w_n)}> n, w_n\ \mbox{not good}\bigg)\nonumber\\
	&\leq \sum_{k=1}^J  \Q\bigg(F_n, \min_{1\leq j<\tau_{\zeta_n}^{(w_n)}} V(w_j)> -\Gamma,  \tau_{\zeta_n}^{(2,w_n)}> n,  \xi(w_k) > e^{V(w_{k-1})+B}  \bigg)\nonumber\\
	&\quad + \sum_{k=J+1}^{n-1} \Q\bigg(F_n, \min_{1\leq j<\tau_{\zeta_n}^{(w_n)}} V(w_j)> -\Gamma,  \tau_{\zeta_n}^{(2,w_n)}> n,  \xi(w_k) > e^{x+ V(w_{k-1})-mk/3} , \tau_{\zeta_n}^{(w_n)}>k \bigg)\nonumber\\
	&=: \sum_{k=1}^{J} R_1(k) + \sum_{k=J+1}^{n-1} R_2(k). 
\end{align}

\paragraph{Bounds for $R_2$.}
We treat $R_2$ first. Combining the Markov property at time $k$ and the fact $\{\min_{1\leq j<\tau_{\zeta_n}^{(w_n)}} V(w_j)> -\Gamma\} \subset \{\min_{1\leq j\leq k} V(w_j)\geq  -\Gamma\}$ on $\{\tau_{\zeta_n}^{(w_n)} >k\}$, it holds that 
\begin{align}
	R_2(k)\leq \E_{\Q}\left(1_{\left\{ \xi(w_k)> e^{x+V(w_{k-1})-mk/3} , \min_{1\leq j\leq k} V(w_j)\geq -\Gamma, \tau_{\zeta_n}^{(w_n)}>k  \right\}} g_{n-k}(V(w_k))\right),
\end{align}
where 
\begin{align}
	& g_{n-k}(z)\nonumber\\
	&:= \mathbf{P}\Big(S_{n-k}\leq \alpha_n -x-z, \min_{\tau_{\zeta_n}\leq j\leq n-k} S_j \geq y_n-z, \tau_{\zeta_n}\in [n-k-T,n-k], \tau_{\zeta_n}^{(2)}> n-k\Big).
\end{align}
Let $n$ be large enough such that $n^{(3-2b)/4}>J+1$ and that $y_n+\Gamma\leq \alpha_n +\Gamma< 4(n-k)^{(3-2b)/4}$ for all $k< n^{(3-2b)/4}$. When $k< n^{(3-2b)/4}$, noticing that on the event $V(w_k)\geq -\Gamma$, we have $y_n-V(w_k)\leq 4(n-k)^{(3-2b)/4}$. Therefore, it follows from Lemma \ref{lem3} (ii) that
\begin{align}\label{Step40}
	g_{n-k}(V(w_k))\lesssim \sum_{z=0}^{[L]} \sum_{i=n-k-T}^{n-k} \left(\mathbf{P}\left(\max_{1\leq j\leq n-k-i} S_j \leq z+1\right) e^{-\alpha_{n-k}} +\frac{1}{n^2}e^{-\alpha_n}\right)\lesssim_{L,T} e^{-\alpha_n},\quad 
\end{align}
where in the last inequality we used the fact that $\alpha_{n-k}=\alpha_n +o(1)$ uniformly on $1\leq k\leq n^{(3-2b)/4}.$ When $n^{(3-2b)/4}\leq k\leq n-1$, on one hand, if $|V(w_k)-mk|\leq n^{(3-2b)/4}$, then 
\begin{align}
	\left|\alpha_n -x- V(w_k)+ mk\right|\leq 2n^{(3-2b)/4}+x\leq 3n^{(3-2b)/4}.
\end{align}
Therefore, by Lemma \ref{lem4}, on the event  $\{|V(w_k)-mk|\leq n^{(3-2b)/4}\}$, it holds that 
\begin{align}\label{Step41}
	g_{n-k}(V(w_k))\lesssim_{L, T} e^{-\alpha_n}.
\end{align}
On the other hand, if $|V(w_k)-mk|> n^{(3-2b)/4}$, then we trivally have $g_{n-k}(z)\leq (1+T)\mathbf{P}(X<-\zeta_n)\lesssim_T \mathbf{P}(X<-\zeta_n)$.
Combining \eqref{Step40}, \eqref{Step41} and the above inequality, we get that 
\begin{align}\label{Step42}
	\sum_{k=J+1}^{n-1} R_2(k)&\lesssim_{L,T} e^{-\alpha_n} \sum_{k=J+1}^{n-1} \Q\Big(\xi(w_k)> e^{x+V(w_{k-1})-mk/3} , \min_{1\leq j\leq k} V(w_j)\geq -\Gamma, \tau_{\zeta_n}^{(w_n)}>k  \Big)\nonumber\\
	&\quad + \mathbf{P}(X<-\zeta_n)\sum_{k=[n^{(3-2b)/4}]+1}^{n-1}  \Q(|V(w_k)-mk|>n^{(3-2b)/4}).
\end{align}
For the last term on the right hand side of \eqref{Step42}, similar to \eqref{Step38}, combining Markov's inequality and  \eqref{Step16}, fixing any integer $D$ such that $(3-2b)D/4 > D/2 +2-b+ \lambda b m^{b-1} A_1 +1$, it holds that 
\begin{align}\label{Step43}
	&\mathbf{P}(X<-\zeta_n) \sum_{k=[n^{(3-2b)/4}]+1}^{n-1} \Q(|V(w_k)-mk|>n^{(3-2b)/4})\nonumber\\
	&\lesssim \sum_{k=[n^{(3-2b)/4}]+1}^{n-T}  e^{-\alpha_n} n^{1-b+ \lambda b m^{b-1} A_1} \frac{k^{D/2}}{n^{(3-2b)D/4}}\nonumber\\
	&\leq  e^{-\alpha_n} n^{2-b+ \lambda b m^{b-1} A_1} \frac{n^{D/2}}{n^{(3-2b)D/4}}\leq \frac{1}{n}e^{-\alpha_n}. 
\end{align}
By the Markov property, $\xi(w_k)$ is independent of $V(w_{k-1})$. Therefore, the first term on the right hand side of \eqref{Step42} has upper bound
\begin{align}\label{Step44}
	&e^{-\alpha_n} \sum_{k=J+1}^{n-1} \Q\left(\xi(w_k)> e^{x+V(w_{k-1})-mk/3} , \min_{j\leq k} V(w_j)\geq -\Gamma, \tau_{\zeta_n}^{(w_n)}>k  \right)\nonumber\\
	&\leq e^{-\alpha_n} \sum_{k=J+1}^{\infty} \Q \left(\log_+ (\xi(w_k))> V(w_{k-1})-mk/3 \right) \nonumber\\
	& = e^{-\alpha_n} \E_{\Q}\left(\sum_{k=J}^{\infty} \mathbf{P}\left( S_{k}-m(k+1)/3 < z\right)\Big|_{z= \log_+(\xi(w_k))} \right).
\end{align}
On the set $ \log_+(\xi(w_k))> mk/6$, we use the trivial upper bound $1$; on the set $ \log_+(\xi(w_k)) \leq mk/6$, by Markov's inequality, when $k$ is large enough such that $m(k+1)/3+ mk/6 \leq 2mk/3$,  we have
\[
\mathbf{P}\left( S_{k}-m(k+1)/3 < z\right) \leq \mathbf{P}\left( |S_{k}-mk| > \frac{mk}{3}\right)\lesssim \frac{1}{k^2}. 
\]
Plugging the above inequality back to \eqref{Step44}, we conclude that when $J$ is large enough,
\begin{align}\label{Step45}
	&e^{-\alpha_n} \sum_{k=J+1}^{n-1} \Q\left(\xi(w_k)> e^{x+V(w_{k-1})-mk/3} , \min_{j\leq k} V(w_j)\geq -\Gamma, \tau_{\zeta_n}^{(w_n)}>k  \right)\nonumber\\
	&\lesssim e^{-\alpha_n} \E_{\Q}\left(\sum_{k=J}^{\infty} 1_{\left\{  \log_+(\xi(w_k)) > mk/6\right\}}\right) + e^{-\alpha_n}\sum_{k=J}^{\infty}  \frac{1}{k^2}\nonumber\\
	&\lesssim e^{-\alpha_n}\left( \E_{\Q}\left(\max\left\{ \log_+ (\xi(w_1))- \frac{6J}{m}, 0 \right\} \right) + \frac{1}{J} \right). 
\end{align}
Therefore, combining \eqref{Step42}, \eqref{Step43} and \eqref{Step45}, we conclude that
\begin{align}\label{Upper-of-R-2}
	\sum_{k=J+1}^{n-1} R_2(k)&\lesssim_{L,T}  e^{-\alpha_n}\left( \E_{\Q}\left(\max\left\{ \log_+ (\xi(w_1))- \frac{6J}{m}, 0 \right\} \right) + \frac{1}{J}+\frac{T+1}{n} \right).  
\end{align}

\paragraph{Bounds for $R_1$.}
Now we deal with $R_1$. Similarly, by the  Markov property at time $k$, for any $1\leq k\leq J$,
\begin{align}
	R_1(k)=\E_{\Q}\left(1_{\left\{ \xi(w_k)> e^{B+V(w_{k-1})} , \min_{1\leq j\leq k} V(w_j)\geq -\Gamma, \tau_{\zeta_n}^{(w_n)}>k  \right\}} g_{n-k}(V(w_k))\right).
\end{align}
Noticing that in this case, when $n$ is large enough, on the set $V(w_k)\geq -\Gamma$, we have
$y_n-V(w_k)\leq \alpha_n+\Gamma <2(n-k)^{(3-2b)/4}$ for all $1\leq k\leq J$. Therefore, combining Lemma \ref{lem3} (ii) and the fact that $e^{\alpha_n -\alpha_{n-k}}\lesssim_J1$ for all $1\leq k\leq J$, we obtain that
\begin{align}\label{Upper-bound-of-R-1}
	\sum_{k=1}^J R_1(k)& \lesssim_J (L+1)e^{-\alpha_n} \sum_{k=1}^J\Q\Big( \xi(w_k)> e^{B+V(w_{k-1})} , \min_{1\leq j\leq k} V(w_j)\geq -\Gamma, \tau_{\zeta_n}^{(w_n)}>k \Big)\nonumber\\
	&\leq (L+1)e^{-\alpha_n} \sum_{k=1}^J \mathbf{P}(S_{k-1}\leq -B/2)+ \frac{(L+1)J}{B}e^{-\alpha_n} \E_{\Q}\left(\log_+ (\xi(w_1))\right). 
\end{align}
Noticing that $\xi(w_1)\leq W_1$, by \eqref{Assumption2}, we have 
\[
\E_{\Q}\left( \log_+ (\xi(w_1))\right)\leq \E_{\Q}\left( \log_+ (W_1)\right)= \E\left(W_1 \log_+ (W_1)\right)<\infty. 
\]
Therefore, for any $\varepsilon, L, T>0$, from \eqref{Upper-of-R-2}, we may take a sufficient large $J$ such that $\sum_{k=J+1}^{n-1} R_2(k)< \frac{\varepsilon}{4}e^{-\alpha_n}$ for large $n$. Then for each given $J$, by \eqref{Upper-bound-of-R-1}, we may take a sufficient large constant $B$ such that $\sum_{k=1}^J R_1(k)<\frac{\varepsilon}{4}e^{-\alpha_n}$, which implies \eqref{Proof-part3}. We are done. 

\hfill$\Box$

We use $\mathbb{T}^{(v)}$ to denote the subtree of $\mathbb{T}$ rooted at $v$. Recall that $R_f>0$ is a fixed constant such that $\mbox{supp}(f)\subset (-\infty, R_f)$. Define
\begin{align}\label{Def-of-E-n-x}
	\mathcal{E}_n(x):=\Big\{\forall k\leq \tau_{\zeta_n}^{(w_n)}, \forall v\in\mathbb{B}(w_k), \min_{u\in\mathbb{T}^{(v)}, |u|=n} V(u)>\alpha_n -x +|R_f| \Big\}.
\end{align}

\begin{prop}\label{prop:bad-event}
	For any $\varepsilon, L>0$,  $T\in \N$ and $R_f\in \R$, there exists $x_*>0$ such that for large $n$, when $x\in [x_*, n^{(3-2b)/4}]$, we have
	\begin{align}
		& \Q\bigg(V(w_n)\leq \alpha_n-x, \min_{\tau_{\zeta_n}^{(w_n)}\leq j\leq n} V(w_j)\geq y_n, \tau_{\zeta_n}^{(w_n)}\in [n-T,n], \left(\mathcal{E}_n(x)\right)^c\bigg)\nonumber\\
		&\leq \varepsilon e^{-\alpha_n}. 
	\end{align}
	Consequently, it holds that 
	\begin{align}
		& \E_{\Q}\bigg( e^{V(w_n)}1_{\left\{V(w_n)\leq \alpha_n-x, \min_{\tau_{\zeta_n}^{(w_n)}\leq j\leq n} V(w_j)\geq y_n, \tau_{\zeta_n}^{(w_n)}\in [n-T,n], \left(\mathcal{E}_n(x)\right)^c \right\}}\bigg)\nonumber\\
		&\leq \varepsilon e^{-x}. 
	\end{align}
\end{prop}
\textbf{Proof: } 
\paragraph{Step 1.}
Recall the definition of $F_n$ in \eqref{Def-of-F-n}. By Lemma \ref{lem:bad-even4}, there exist two large constants $J,B$ such that 
\begin{align}
	\Q\left(F_n, \left(\mathcal{E}_n(x)\right)^c\right) & \leq  \frac{\varepsilon}{4} e^{-\alpha_n}+ \Q\left(F_n, \left(\mathcal{E}_n(x)\right)^c, w_n\ \mbox{good}\right). 
\end{align}
Therefore, it mains to show that when $n$ is large enough,
\begin{align}\label{Goal2}
	\Q\left(F_n, \left(\mathcal{E}_n(x)\right)^c, w_n\ \mbox{good}\right)\leq \frac{3\varepsilon}{4}e^{-\alpha_n}. 
\end{align}
Noticing that when $n$ is large enough such that $n^{(3-2b)/4}/2 > mT$, we have
\begin{align}
	& \sum_{k=n-T}^n \Q\left(F_n, \tau_{\zeta_n}^{(w_n)}=k, \tau_{\zeta_n}^{(2,w_n)}>n, |V(w_{k-1})-m(n-1)|>n^{(3-2b)/4} \right)\nonumber\\
	&\leq \sum_{k=n-T}^n \mathbf{P}(X<-\zeta_n)\mathbf{P}\left( |S_{k-1}-m(k-1)|>\tfrac{1}{2}n^{(3-2b)/4} \right).
\end{align}
Therefore, fixing any integer $D$ such that $(3-2b)D/4> D/2 + 2-b+\lambda bm^{b-1} A_1$, by \eqref{Step16},
\begin{align}\label{Step46}
	& \sum_{k=n-T}^n \Q\left(F_n, \tau_{\zeta_n}^{(w_n)}=k, \tau_{\zeta_n}^{(2,w_n)}>n, |V(w_{k-1})-m(n-1)|\geq n^{(3-2b)/4} \right)\nonumber\\
	&\lesssim (T+1)e^{-\alpha_n} n^{1-b+ \lambda b m^{b-1} A_1} \frac{n^{D/2}}{n^{(3-2b)D/4}}\leq \frac{T+1}{n}e^{-\alpha_n}. 
\end{align}
According to Lemma \ref{lem:bad-even5}, for large $n$, 
\begin{align}\label{Step47}
	& \sum_{k=n-T}^n \Q\left(F_n, \tau_{\zeta_n}^{(w_n)}=k, \tau_{\zeta_n}^{(2,w_n)}>n, |V(w_k)-(\alpha_n -x)|>\log n \right)\nonumber\\
	&\leq \sum_{z=0}^{[L]} \mathbf{P}\left(S_n -y_n\in [z, z+1], \min_{\tau_{\zeta_n}  \leq j\leq n} S_j\geq y_n, \tau_{\zeta_n} \leq n, |S_n- S_{\tau_{\zeta_n} }|>\log n-L \right)\nonumber\\
	&\leq (L+1) \sup_{z\in [0,L]} \mathbf{P}\left(S_n -y_n\in [z, z+1], \min_{\tau_{\zeta_n}  \leq j\leq n} S_j\geq y_n, \tau_{\zeta_n} \leq n, |S_n- S_{\tau_{\zeta_n} }|>\log n-L \right)\nonumber\\
	&\leq \frac{\varepsilon}{8}e^{-\alpha_n}.
\end{align}
Now define 
\begin{align}
	\Upsilon_k^{(n)}:= \left\{|V(w_{k-1})-m(n-1)|<n^{(3-2b)/4}, \ |V(w_k)-(\alpha_n -x)| \leq \log n \right\}.
\end{align}
Then combining \eqref{Step46} and \eqref{Step47}, we see that when $n$ is large enough,
\begin{align}\label{Step48}
	\sum_{k=n-T}^n \Q\left(F_n, \tau_{\zeta_n}^{(w_n)}=k, \tau_{\zeta_n}^{(2,w_n)}>n,  \left(\Upsilon_k^{(n)}\right)^c \right)\leq \frac{\varepsilon}{4}e^{-\alpha_n}. 
\end{align}
According to the definition of $\mathcal{E}_n(x)$, we have
\begin{align}\label{Proof-Goal2-1}
	&\Q\left(F_n, \left(\mathcal{E}_n(x)\right)^c, w_n\ \mbox{good}\right)\nonumber\\
	&\leq \sum_{k=n-T}^n \Q\left(F_n, \tau_{\zeta_n}^{(w_n)}=k, \tau_{\zeta_n}^{(2,w_n)}>n,  \exists v\in \mathbb{B}(w_k), \min_{u\in\mathbb{T}^{(v)}, |u|=n} V(u)\leq \alpha_n-x +|R_f| \right)\nonumber\\
	&\quad+\sum_{k=n-T}^n \E_{\Q}\bigg(1_{\left\{F_n, \tau_{\zeta_n}^{(w_n)}=k, w_n\ \mbox{good} \right\}} \sum_{j=1}^{k-1} \sum_{v\in\mathbb{B}(w_j)} 1_{\left\{\min_{u\in\mathbb{T}^{(v)}, |u|=n} V(u)\leq \alpha_n -x +|R_f| \right\}}\bigg)\nonumber\\
	&=: U_n+V_n.
\end{align}

\paragraph{Step 2.} In this step, we prove an upper bound for $U_n$.  Combining \eqref{Step48} and the Markov property,  when $n$ is large enough such that $|R_f|\leq n^{(3-2b)/4}$, 
\begin{align}\label{Step52}
	U_n&\leq \frac{\varepsilon}{4}e^{-\alpha_n}+ \sum_{k=n-T}^n \Q\left(\Upsilon_k^{(n)}, \tau_{\zeta_n}^{(w_n)}=k,  \exists v\in \mathbb{B}(w_k), \min_{u\in\mathbb{T}^{(v)}, |u|=n} V(u)\leq \alpha_n-x+|R_f| \right)\nonumber\\
	&\leq \frac{\varepsilon}{4}e^{-\alpha_n}+ \sum_{k=0}^T \sup_{|z+m(n-1)|<4n^{(3-2b)/4}}\Q\bigg(|V(w_1)-z|\leq \log n+|R_f|, \nonumber\\
	&\qquad\qquad\qquad \exists v\in \mathbb{B}(w_1), \min_{u\in\mathbb{T}^{(v)}, |u|=k+1} V(u)\leq -(mn-5n^{(3-2b)/4})\bigg),
\end{align}
where in the last inequality we used the fact that for  $z:=- V(w_{k-1})+\alpha_n -x+|R_f| $, we have $|z+m(n-1)|<4n^{(3-2b)/4}$ and that $\alpha_n -x- V(w_{k-1}) +|R_f|\leq -(mn-5n^{(3-2b)/4})$ on the set $\Upsilon_k^{(n)}$ when $n$ is large enough. According to the many-to-one formula \eqref{Many-to-one} and the definition of $\mathscr{L}$, we see that 
\begin{align}\label{Step49}
	&\Q\left(|V(w_1)-z|\leq \log n+|R_f|, \exists v\in \mathbb{B}(w_1), \min_{u\in\mathbb{T}^{(v)}, |u|=k+1} V(u)\leq -(mn-5n^{(3-2b)/4})\right)\nonumber\\
	& = \E\left(\sum_{i=1}^\nu e^{-Y_i}1_{\left\{|Y_i-z|\leq \log n+|R_f| \right\}} 1_{\left\{\exists j\neq i, Y_j +M_k^{(j)}\leq -(mn-5n^{(3-2b)/4})\right\}}\right)\nonumber\\
	&= \sum_{\ell=1}^\infty \ell \P(\nu=\ell) \E\left( e^{-Y_1}1_{\left\{|Y_1-z|\leq \log n + |R_f| \right\}} 1_{\left\{ \exists j\in \{2,...,\ell \}, Y_j +M_k^{(j)}\leq -(mn-5n^{(3-2b)/4})\right\}}\right)\nonumber\\
	&\leq \sum_{\ell=1}^\infty \ell \P(\nu=\ell) \E\Big( e^{-Y_1}1_{\left\{|Y_1-z|\leq \log n+|R_f| \right\}}\Big) \nonumber\\
	&\quad\quad\times \P\left(  \exists j\in \{2,...,\ell \}, Y_j +M_k^{(j)}\leq -(mn-5n^{(3-2b)/4})\right) , 
\end{align}
where $M_k^{(j)}$ are iid equal in law to $M_k$ under $\P$ and in the last inequality we also used the fact that $Y_i$ are independent to each other. Combining the independence of $Y_i$ and $\nu$ and many-to-one formula \eqref{Many-to-one}, 
\begin{align}\label{Step50}
	& \E\left( e^{-Y_1}1_{\left\{|Y_1-z|\leq \log n +|R_f|\right\}}\right)  = \frac{1}{\E \nu} \E\left(\sum_{i=1}^\nu  e^{-Y_i}1_{\left\{|Y_i-z|\leq \log n +|R_f| \right\}}\right) \nonumber\\
	&= \frac{1}{\E \nu} \mathbf{P}(|X-z|\leq \log n+|R_f| )\nonumber\\
	&\lesssim \int_{z-\log n-|R_f|}^{z+\log n+|R_f| }|y|^a e^{-\lambda |y|^b}\mathrm{d} y \lesssim( \log n)e^{-\alpha_n},
\end{align}
here in the last inequality we used the fact that $|z+m(n-1)|< 4n^{(3-2b)/4}$.  Moreover, if $\ell \leq e^{mn/2}$, then by \eqref{Tail-lowest}, we see that 
\begin{align}\label{Step51}
	& \P\left(  \exists j\in \{2,...,\ell \}, Y_j +M_k^{(j)} \leq -(mn-5n^{(3-2b)/4})\right)\nonumber\\
	&\leq \sum_{j=2}^\ell  \E\left(e^{-Y_j-(mn-5n^{(3-2b)/4} )} \right)\leq \E(e^{-Y_1}) e^{-(mn/2-5n^{(3-2b)/4} )}\lesssim e^{-(mn/2-5n^{(3-2b)/4})}.
\end{align}
Now combining \eqref{Step49}, \eqref{Step50} and \eqref{Step51}, we obtain that
\begin{align}\label{Step53}
	&\Q\left(|V(w_1)-z|\leq \log n +|R_f|, \exists v\in \mathbb{B}(w_1), \min_{u\in\mathbb{T}^{(v)}, |u|=k+1} V(u)\leq -(mn-5n^{(3-2b)/4})\right)\nonumber\\
	&\lesssim (\log n) e^{-\alpha_n} \bigg( e^{-(mn/2-5n^{(3-2b)/4})} \sum_{\ell=1}^{[e^{mn/2}]} \ell \P(\nu=\ell) +\sum_{\ell=[e^{mn/2}]}^{\infty} \ell \P(\nu=\ell) \bigg)\nonumber\\
	&\lesssim \frac{\log n}{n} e^{-\alpha_n} \E (\nu \log_+ \nu ). 
\end{align}
Since $\E (\nu \log_+ \nu )<\infty$ by \eqref{Assumption2} (or see \eqref{Equal-Conditions}), combining \eqref{Step52} and \eqref{Step53}, when $n$ is large enough, we have
\begin{align}\label{Proof-Goal2-2}
	U_n \leq \frac{\varepsilon}{2}e^{-\alpha_n}. 
\end{align}

{\bf(Step 3).}  In this step, we treat $V_n$. Let $\mathcal{G}$ be the $\sigma$-field generated by all the information along the spine and their siblings. Then by the Markov property, for any $v\in \mathbb{B}(w_j)$, 
\begin{align}\label{Step54}
	\Q\left(\min_{u\in\mathbb{T}^{(v)}, |u|=n} V(u)\leq \alpha_n -x +|R_f| \bigg|\mathcal{G} \right) = \P(M_{n-j}\leq \alpha_n-x-z+|R_f|)\big|_{z=V(v)}.
\end{align}
If $j\leq n^{(3-2b)/4}$, then by Proposition \ref{prop:Tightness}, and the fact that $\alpha_{n-j}=\alpha_n+o(1)$, we have
\begin{align}\label{Step55}
	\P(M_{n-j}\leq \alpha_n-x-z+|R_f|)\lesssim_{R_f} e^{-(z+x+\alpha_{n-j}-\alpha_n)}\lesssim  e^{-(z+x)}. 
\end{align}
If $j>n^{(3-2b)/4}$, then by \eqref{Tail-lowest}, we get that 
\begin{align}\label{Step56}
	\P(M_{n-j}\leq \alpha_n-x-z+|R_f|)\lesssim_{R_f} e^{-(z+x-\alpha_n)}.
\end{align}
Combining \eqref{Step54}, \eqref{Step55} and \eqref{Step56}, we get that when $n$ is large enough, on the set $\{w_n \mbox{ is good}\}$, it holds that 
\begin{align}\label{Step57}
	& \E_{\Q}\bigg( \sum_{v\in\mathbb{B}(w_j)} 1_{\left\{\min_{u\in\mathbb{T}^{(v)}, |u|=n} V(u)\leq \alpha_n -x +|R_f| \right\}}\bigg|\mathcal{G}\bigg) \nonumber\\
	&\lesssim_{R_f} \begin{cases}
		e^{B-x},\quad & \mbox{if} \ 1\leq j\leq J;\\
		e^{-mj/3} ,\quad &\mbox{if} \ J<j\leq n^{(3-2b)/4}; \\
		e^{\alpha_n-mj/3},\quad &\mbox{if} \ n^{(3-2b)/4}<j<k. 
	\end{cases}
\end{align}
Combining the definition of $V_n$ in \eqref{Proof-Goal2-1} and  \eqref{Step57},  we see that 
\begin{align}\label{Step58}
	V_n & \lesssim_{R_f} \sum_{k=n-T}^n \E_{\Q}\bigg(1_{\left\{F_n, \tau_{\zeta_n}^{(w_n)}=k, w_n\ \mbox{good} \right\}} \Big(\sum_{j=1}^{J} e^{B-x}+ \sum_{j=J+1}^{[n^{(3-2b)/4}]} e^{-mj/3} + \sum_{j=[n^{(3-2b)/4}]+1}^{k-1} 	e^{\alpha_n-mj/3} \Big)\bigg)\nonumber\\
	&\leq \sum_{k=n-T}^n \Q\left(F_n, \tau_{\zeta_n}^{(w_n)}=k\right)  \bigg(J e^{B-x}+ \sum_{j=J+1}^{\infty} e^{-mj/3} + \sum_{j=[n^{(3-2b)/4}]+1}^{n} 	e^{-mj/6} \bigg),
\end{align}
where in the last inequality we used the fact that $\alpha_n \leq m n^{(3-2b)/4} /6$ when $n$ is large enough.  By Lemma \ref{lem3} (ii), we obtain
\[
\sum_{k=n-T}^n \Q\left(F_n, \tau_{\zeta_n}^{(w_n)}=k\right) \lesssim (T+1)(L+1)e^{-\alpha_n},
\]
which together with \eqref{Step58} implies that 
\begin{align}
	V_n \lesssim_{R_f}  \bigg(J e^{B-x}+ \sum_{j=J+1}^{\infty} e^{-mj/3} + \sum_{j=[n^{(3-2b)/4}]+1}^{n} 	e^{-mj/6} \bigg) (T+1)(L+1)e^{-\alpha_n}. 
\end{align}
Therefore, taking $J$ sufficient large first and then $x_*$ sufficient large, for large $n$ and $x\geq x_*$, 
\begin{align}\label{Proof-Goal2-3}
	V_n \leq \frac{\varepsilon}{4}e^{-\alpha_n}. 
\end{align}
Now \eqref{Goal2} follows directly from \eqref{Proof-Goal2-1}, \eqref{Proof-Goal2-2} and \eqref{Proof-Goal2-3}, which implies the desired result. 

\hfill$\Box$

\textbf{Proof of Proposition \ref{prop:key-proposition}: }
Noticing that 
\begin{align}\label{Step74}
	& \E \left(1- e^{-\sum_{|u|=n} f(V(u)- (\alpha_n -x))} 1_{\left\{ M_n> \alpha_n -x \right\}}\right)\nonumber\\
	& = \E \left(1- e^{-\sum_{|u|=n} f(V(u)- (\alpha_n -x))} \right) + \E \left( e^{-\sum_{|u|=n} f(V(u)- (\alpha_n -x))} 1_{\left\{ M_n\leq \alpha_n -x \right\}}\right)\nonumber\\
	& = \E \left(\left(1- e^{-\sum_{|u|=n} f(V(u)- (\alpha_n -x))}\right)1_{\left\{M_n\leq \alpha_n -x+R_f \right\}} \right) + \E \left( e^{-\sum_{|u|=n} f(V(u)- (\alpha_n -x))} 1_{\left\{ M_n\leq \alpha_n -x \right\}}\right)\nonumber\\
	&=: \E \left(Z^{(1)}_{n,x_1} 1_{\left\{M_n\leq \alpha_n -x_1 \right\}} \right) + \E \left( Z^{(2)}_{n,x_2} 1_{\left\{ M_n\leq \alpha_n -x_2 \right\}}\right),
\end{align}
where 
\begin{align}
	x_1& := x-R_f, \quad Z^{(1)}= Z^{(1)}_{n,x_1} := 1-e^{-\sum_{|u|=n} f(V(u)- (\alpha_n -x_1-R_f))} ,\nonumber\\
	&x_2=x\quad \mbox{and}\quad Z^{(2)}= Z^{(2)}_{n,x_2} := e^{-\sum_{|u|=n} f(V(u)- (\alpha_n -x_2))}. 
\end{align}
We will treat both of the cases simultaneously. 
For any $\varepsilon_0>0$, by Lemma \ref{lem:bad-event1}, when $n$ is large enough, for any $x_i \geq 0$, $i=1,2$, 
\begin{align}
	\bigg|	\E \left(Z^{(i)}1_{\left\{ M_n\leq  \alpha_n -x_i \right\}}\right) - 	\E \left(Z^{(i)}1_{\left\{ \exists u\in \mathbb{T}_n , V(u)=M_n\leq \alpha_n -x_i, \mathcal{H}_{\zeta_n, n}^u \right\}} \right) \bigg|\leq 2\varepsilon_0 e^{-x_i}. 
\end{align}
Also applying \eqref{Eq2} in Lemma \ref{lem:bad-event2}, taking $L_0=L_0(\varepsilon)$ sufficiently large, we obtain from the above inequality that for any $L\geq L_0$, when $n$ is large enough, for any $x_i\geq 0$, $i=1,2$,
\begin{align}\label{Step59}
	& \bigg|	\E \left(Z^{(i)}1_{\left\{ M_n\leq \alpha_n -x_i \right\}}\right)  \nonumber\\
	&\quad - 	\E\bigg(Z^{(i)}1_{\left\{ \exists u\in \mathbb{T}_n , V(u)=M_n\leq \alpha_n -x_i, \min_{\tau_{\zeta_n}^{(u)}\leq j\leq n} V(u_j) >\alpha_n-x_i-L,  \mathcal{H}_{\zeta_n, n}^u \right\}} \bigg) \bigg| \leq 3\varepsilon_0 e^{-x_i}. 
\end{align}
Combining \eqref{Step59} and Lemma \ref{lem:bad-event3}, we see that for any $L\geq L_0$, there exists an integer $T_0=T_0(\varepsilon_0, L)$ such that for any $T\geq T_0$, when $n$ is large enough, for all $x_i\geq 0$, $i=1,2$, 
\begin{align}\label{Step60}
	& \bigg|\E \left(Z^{(i)}1_{\left\{ M_n\leq \alpha_n -x_i \right\}}\right)  \nonumber\\
	&\quad - 	\E \bigg(Z^{(i)}1_{\left\{ \exists u\in \mathbb{T}_n , V(u)=M_n\leq \alpha_n -x_i, \min_{\tau_{\zeta_n}^{(u)}\leq j\leq n} V(u_j) >\alpha_n-x_i-L,   n-T\leq \mathcal{H}_{\zeta_n, n}^u \right\}}\bigg) \bigg| \leq 4\varepsilon_0 e^{-x_i}. 
\end{align}
Now again by \eqref{Eq-at-most-two-jumps-2} in  Lemma \ref{lem:bad-event1}, we can drop out the event $\{\tau_{\zeta_n}^{(2,u)}>n\}$ in \eqref{Step60} to get that 
\begin{align}\label{Step63}
	& \bigg|	\E \left(Z^{(i)}1_{\left\{ M_n\leq \alpha_n -x_i \right\}}\right) \nonumber\\
	&\quad - 	\E \bigg(Z^{(i)}1_{\left\{ \exists u\in \mathbb{T}_n , V(u)=M_n\leq \alpha_n -x_i, \min_{\tau_{\zeta_n}^{(u)}\leq j\leq n} V(u_j) >\alpha_n-x_i-L,   n-T\leq \tau_{\zeta_n}^{(u)}\leq  n\right\}}\bigg) \bigg| \leq 5\varepsilon_0 e^{-x_i}. 
\end{align}
Combining \eqref{Change-of-measure} and \eqref{Gibbs-measure}, 
\begin{align}\label{Step73}
	&	\E \bigg(Z^{(i)}1_{\left\{ \exists u\in \mathbb{T}_n , V(u)=M_n\leq \alpha_n -x_i, \min_{\tau_{\zeta_n}^{(u)}\leq j\leq n} V(u_j) >\alpha_n-x_i-L,   n-T\leq \tau_{\zeta_n}^{(u)}\leq  n\right\}}\bigg)  \nonumber\\
	& = \E\bigg( \frac{Z^{(i)}}{\sum_{|u|=n}1_{\{V(u)=M_n\}} } \sum_{|u|=n}  1_{\left\{V(u)=M_n\leq 
		\alpha_n-x_i, \min_{\tau_{\zeta_n}^{(u)}\leq j\leq n} V(u_j) >\alpha_n-x_i-L,   n-T\leq \tau_{\zeta_n}^{(u)}\leq  n \right\}}\bigg)\nonumber\\
	& = \E_{\Q} \bigg( e^{V(w_n)}Z^{(i)} \frac{ 1_{\{V(w_n)=M_n\leq 
			\alpha_n-x_i, \min_{\tau_{\zeta_n}^{(w_n)}\leq j\leq n} V(w_j) >\alpha_n-x_i-L,   n-T\leq \tau_{\zeta_n}^{(w_n)}\leq  n\}}  }{\sum_{|u|=n}1_{\{V(u)=M_n\}} }  \bigg).
\end{align}
Recall the definition of $\mathcal{E}_n(x)$ in \eqref{Def-of-E-n-x}. On the set $\mathcal{E}_n(x)$, we know that if $k=\tau_{\zeta_n}^{(w_n)}$, then
\begin{align}
	& Z^{(1)} = Z_{n,x_1}^{(1)}(w_k):= 1-e^{-\sum_{u\in \mathbb{T}^{(w_k)}, |u| =n} f(V(u)-(\alpha_n-x_1-R_f))} ,\nonumber\\
	& Z^{(2)}=Z_{n,x_2}^{(2)}(w_k):= e^{-\sum_{u\in \mathbb{T}^{(w_k)}, |u| =n} f(V(u)-(\alpha_n-x_2))} ,\nonumber\\
	& \sum_{|u|=n}1_{\{V(u)=M_n\}} =\sum_{u\in \mathbb{T}^{(w_k)}, |u|=n}1_{\{V(u)=M_n\}} . 
\end{align}
Therefore, combining  \eqref{Step63}, \eqref{Step73} and Proposition \ref{prop:bad-event}, there exists $x_*>0$ such that for $n$ large enough, for any $x_i\in [x_*, n^{(3-2b)/4}]$, $i=1,2$, 
\begin{align}\label{Step61}
	& \bigg|	\E \left(Z^{(i)}1_{\left\{ M_n\leq \alpha_n -x_i \right\}}\right)  \nonumber\\
	&\quad- 	\E_{\Q} \bigg( e^{V(w_n)} Z_{n,x_i}^{(i)}(w_k) \frac{ 1_{\{V(w_n)=M_n\leq 
			\alpha_n-x_i, \min_{\tau_{\zeta_n}^{(w_n)}\leq j\leq n} V(w_j) >\alpha_n-x_i-L,   n-T\leq \tau_{\zeta_n}^{(w_n)}\leq  n, \mathcal{E}_n(x) \}}  }{\sum_{u\in \mathbb{T}^{(w_k)}, |u|=n}1_{\{V(u)=M_n\}} }  \bigg) \bigg|\nonumber\\
	&\leq 6\varepsilon_0 e^{-x_i}.
\end{align}
Again by Proposition \ref{prop:bad-event}, we can drop out the event $\mathcal{E}_n(x)$ in \eqref{Step61} and get that 
\begin{align}\label{Step62}
	& \bigg|	\E \left(Z^{(i)}1_{\left\{ M_n\leq \alpha_n -x_i \right\}}\right)   \nonumber\\
	&\quad- 	\E_{\Q} \bigg( e^{V(w_n)}Z_{n,x_i}^{(i)}(w_k) \frac{ 1_{\{V(w_n)=M_n\leq 
			\alpha_n-x_i, \min_{\tau_{\zeta_n}^{(w_n)}\leq j\leq n} V(w_j) >\alpha_n-x_i-L,   n-T\leq \tau_{\zeta_n}^{(w_n)}\leq  n \}}  }{\sum_{u\in \mathbb{T}^{(w_k)}, |u|=n}1_{\{V(u)=M_n\}} }  \bigg) \bigg|\nonumber\\
	&\leq 7\varepsilon_0 e^{-x_i}. 
\end{align}
Since $\{ \tau_{\zeta_n}^{(w_n)}=k \}, $ it holds that $\tau_{\zeta_n}^{(w_n)} =\tau_{\zeta_n}^{(w_k)}$. Therefore,  by the Markov property at $k$, 
\begin{align}\label{Step69}
	&\E_{\Q} \bigg( e^{V(w_n)}Z_{n,x_i}^{(i)}(w_k) \frac{ 1_{\{V(w_n)=M_n\leq 
			\alpha_n-x_i, \min_{\tau_{\zeta_n}^{(w_n)}\leq j\leq n} V(w_j) >\alpha_n-x_i-L,   n-T\leq \tau_{\zeta_n}^{(w_n)}\leq  n \}}  }{\sum_{u\in \mathbb{T}^{(w_k)}, |u|=n}1_{\{V(u)=M_n\}} }  \bigg) \nonumber\\
	&= \sum_{k=n-T}^n \E_{\Q} \bigg( e^{V(w_n)}Z_{n,x_i}^{(i)}(w_k)  \frac{ 1_{\{V(w_n)=M_n\leq 
			\alpha_n-x_i, \min_{k\leq j\leq n} V(w_j) >\alpha_n-x_i-L,    \tau_{\zeta_n}^{(w_n)}=k \}}  }{\sum_{u\in \mathbb{T}^{(w_k)}, |u|=n}1_{\{V(u)=M_n\}} }  \bigg)\nonumber\\
	&= e^{\alpha_n -x_i}\sum_{k=n-T}^n \E_{\Q} \bigg(1_{\left\{\tau_{\zeta_n}^{(w_k)}=k, V(w_k)\geq \alpha_n-x_i-L \right\}} F_{n-k}^{(i, L)}\left(V(w_k)- (\alpha_n-x_i-L)\right) \bigg)\nonumber\\
	& = e^{\alpha_n -x_i}\sum_{k=0}^T \E_{\Q} \bigg(1_{\left\{\tau_{\zeta_n}^{(w_{n-k})}=n-k, V(w_{n-k})\geq \alpha_n-x_i-L \right\}} F_{k}^{(i, L)}\left(V(w_{n-k})- (\alpha_n-x_i-L)\right) \bigg), \quad 
\end{align}
where 
\begin{align}
	F_j^{(1,L)}(s)&:=  e^{s-L}\E_{\Q}\left(e^{V(w_j)} \left(1-e^{-\sum_{|v|=j} f(V(v)+s-L+R_f)}  \right) \frac{1_{\{V(w_j)=M_j\leq 
			L-s, \min_{ \ell \leq j} V(w_\ell) >-s \}}  }{\sum_{|v|=j}1_{\{V(v)=M_j\}} } \right).    \nonumber\\
	F_j^{(2, L)}(s)&:= e^{s-L}\E_{\Q}\left(e^{V(w_j)} e^{-\sum_{|v|=j} f(V(v)+s-L)}   \frac{1_{\{V(w_j)=M_j\leq 
			L-s, \min_{ \ell \leq j} V(w_\ell) >-s \}}  }{\sum_{|v|=j}1_{\{V(v)=M_j\}} } \right).
\end{align}
Noticing that  for any $s\geq 0$ and $i=1,2$, 
\begin{align}
	F_j^{(i, L)}(s)& \leq e^{s-L}\E_{\Q}\left(e^{V(w_j)} 1_{\{V(w_j)\leq 
		L-s, \min_{ \ell \leq j} V(w_\ell) >-s \}}   \right)\nonumber\\
	& \leq \mathbf{P}\left(S_j \leq L-s\right)\leq j \mathbf{P}(X< (L-s)/j ),
\end{align}
since $\mathbf{E}(|X|^k)<\infty$ for any $k\in\N$ by \eqref{Assumption4}, we see that for any $T, L>0$ and $i=1,2$, $F_j^{(i, L)}(s)$ satisfies the condition of Lemma \ref{lem:Renewal-theorem} with $H(\zeta,z)=H_\infty(\zeta)= F_j^{(i, L)}(\zeta)$ for all $1\leq j\leq T$. Therefore, combining Lemma \ref{lem:Renewal-theorem}, \eqref{Step62} and \eqref{Step69}, for any $L\geq L_0$ and sufficiently large $T\geq T_0=T_0(\varepsilon, L),$ when $n$ is large enough, for $x\in [x_*, n^{(3-2b)/4}]$, 
\begin{align}\label{Step70}
	& \bigg|			\E \left(Z^{(i)}1_{\left\{ M_n\leq \alpha_n -x_i \right\}}\right)    -	e^{-x_i}\ell_\infty m^a  \sum_{j=0}^T \int_0^\infty F_j^{(i, L)}(\zeta)\mathrm{d}\zeta \bigg|\leq 8\varepsilon_0 e^{-x_i}. 
\end{align}
Since $F_j^{(i, L)}$ is non-negative, by the Fubini theorem, we see that for $i=1$, 
\begin{align}
	& \int_0^\infty F_j^{(1, L)}(\zeta)\mathrm{d}\zeta \nonumber\\
	& = \E_{\Q} \bigg( e^{V(w_j)}  \frac{1_{\{V(w_j)=M_j\}}  }{\sum_{|v|=j}1_{\{V(v)=M_j\}} } \int_{(-\min_{\ell\leq j} V(w_\ell))\vee 0} ^{L-V(w_j)} e^{\zeta-L} \left(1-e^{-\sum_{|v|=j} f(V(v)+\zeta-L+R_f)}  \right)  \mathrm{d}\zeta\bigg)\nonumber\\
	& =  \E_{\Q} \bigg(   \frac{1_{\{V(w_j)=M_j\}}  }{\sum_{|v|=j}1_{\{V(v)=M_j\}} } \int_{0} ^{L-V(w_j)+(\min_{\ell\leq j} V(w_\ell) )\land 0} e^{-\zeta}\left(1-e^{-\sum_{|v|=j} f(V(v)-V(w_j)-\zeta+R_f)} \right) \mathrm{d}\zeta\bigg)
\end{align}
and for $i=2$, 
\begin{align}
	& \int_0^\infty F_j^{(2, L)}(\zeta)\mathrm{d}\zeta = \E_{\Q} \left( e^{V(w_j)}  \frac{1_{\{V(w_j)=M_j\}}  }{\sum_{|v|=j}1_{\{V(v)=M_j\}} } \int_{(-\min_{\ell\leq j} V(w_\ell))\vee 0} ^{L-V(w_j)} e^{\zeta-L}e^{-\sum_{|v|=j} f(V(v)+\zeta-L)}  \mathrm{d}\zeta\right)\nonumber\\
	& =  \E_{\Q} \left(   \frac{1_{\{V(w_j)=M_j\}}  }{\sum_{|v|=j}1_{\{V(v)=M_j\}} } \int_{0} ^{L-V(w_j)+(\min_{\ell\leq j} V(w_\ell) )\land 0} e^{-\zeta}e^{-\sum_{|v|=j} f(V(v)-V(w_j)-\zeta)}  \mathrm{d}\zeta\right).
\end{align}
Now by Proposition \ref{prop:Tightness},  $\sum_{j=0}^T \int_0^\infty F_j^{(i, L)}(\zeta)\mathrm{d}\zeta $ is uniformly bounded for large $T$ and $L$. Therefore, by monotonicity convergence theorem, we obtain that 
\begin{align}\label{Step71}
	& \lim_{L\to\infty} \lim_{T\to\infty} \sum_{j=0}^T \int_0^\infty F_j^{(1, L)}(\zeta)\mathrm{d}\zeta \nonumber\\
	&=  \sum_{j=0}^\infty \E_{\Q} \left(   \frac{1_{\{V(w_j)=M_j\}}  }{\sum_{|v|=j}1_{\{V(v)=M_j\}} } \int_{0} ^{\infty } e^{-\zeta}\left(1-e^{-\sum_{|v|=j} f(V(v)-V(w_j)-\zeta+R_f)} \right) \mathrm{d}\zeta\right)\nonumber\\
	& =  \sum_{j=0}^\infty \E \left(  e^{-M_j}  \int_{0} ^{\infty} e^{-\zeta}\left(1- e^{-\sum_{|v|=j} f(V(v)-M_j-\zeta+R_f)}  \right)\mathrm{d}\zeta\right),
\end{align}
where in the last inequality we used many-to-one formula \eqref{Many-to-one}, i.e., 
\begin{align}
	& \E \left(  e^{-M_j}  \int_{0} ^{\infty} e^{-\zeta}\left(1- e^{-\sum_{|v|=j} f(V(v)-M_j-\zeta+R_f)}\right)  \mathrm{d}\zeta\right)\nonumber\\
	&  = \E\bigg(\frac{1}{\sum_{|u|=j}1_{\{V(u)=M_j\}} } \sum_{|u|=j} e^{-V(u)}1_{\{V(u)=M_j\}}  \int_{0} ^{\infty} e^{-\zeta}\left(1- e^{-\sum_{|v|=j} f(V(v)-V(u)-\zeta+R_f)} \right) \mathrm{d}\zeta\bigg)\nonumber\\
	& =  \E_{\Q} \left(   \frac{1_{\{V(w_j)=M_j\}}  }{\sum_{|v|=j}1_{\{V(v)=M_j\}} } \int_{0} ^{\infty} e^{-\zeta}\left(1- e^{-\sum_{|v|=j} f(V(v)-V(w_j)-\zeta+R_f)} \right) \mathrm{d}\zeta\right) . 
\end{align}
Therefore, combining \eqref{Step70} and \eqref{Step71}, for each given $\varepsilon_0>0$,  we can choose suitable sufficiently large $T$, $L$ and a suitable $A=x_*(T,L, R_f)$ such that for any large $n$, when $x\in [A, n^{(3-2b)/4}],$
\begin{align}\label{Step75}
	& \bigg|	\E \left(Z^{(1)}1_{\left\{ M_n\leq \alpha_n -x_1 \right\}}\right)    -	e^{-x_1}\ell_\infty m^a    \sum_{j=0}^\infty \E \left(  e^{-M_j}  \int_{0} ^{\infty} e^{-\zeta}\left(1- e^{-\sum_{|v|=j} f(V(v)-M_j-\zeta+R_f)}  \right)\mathrm{d}\zeta\right) \bigg|\nonumber\\
	&\leq 9\varepsilon_0 e^{-x_1}.
\end{align} 
Noticing that we have the following the following identity
\begin{align}
	& e^{R_f}\ell_\infty m^a    \sum_{j=0}^\infty \E  \left(  e^{-M_j}  \int_{0} ^{\infty} e^{-\zeta}\left(1- e^{-\sum_{|v|=j} f(V(v)-M_j-\zeta+R_f)}  \right)\mathrm{d}\zeta\right) \nonumber\\
	& = \ell_\infty m^a    \sum_{j=0}^\infty \E  \left(  e^{-M_j}  \int_{\R}  e^{-\zeta}\left(1- e^{-\sum_{|v|=j} f(V(v)-M_j-\zeta)}  \right)\mathrm{d}\zeta\right) = C^*(f)-C^*(0),
\end{align}
where in the first equality we also used the fact that for $\zeta\leq 0$, $V(v)-M_j-\zeta+R_f\geq -\zeta+R_f\geq R_f$ and this implies $1- e^{-\sum_{|v|=j} f(V(v)-M_j-\zeta+R_f)}=0 $.  Therefore, combining the definitions of $x_1$ and $Z^{(1)}$ and \eqref{Step75}, we get that
\begin{align}\label{Step76'}
	& \bigg|	\E\left( \left(1- e^{-\sum_{|u|=n} f(V(u)- (\alpha_n -x))} \right)1_{\{M_n\leq \alpha_n-x_1\}} \right)  -	e^{-x}\left(C^*(f)-C^*(0)\right)  \bigg| \leq 9\varepsilon_0 e^{-x},
\end{align}
which implies \eqref{Second-Ineq}. 
Similarly for $i=2$, we also have
\begin{align}\label{Step76}
	& \bigg|	\E \left(Z^{(2)}1_{\left\{ M_n\leq \alpha_n -x_2 \right\}}\right)    -	e^{-x_2}\ell_\infty m^a    \sum_{j=0}^\infty \E  \left(  e^{-M_j}  \int_{0} ^{\infty} e^{-\zeta}e^{-\sum_{|v|=j} f(V(v)-M_j-\zeta)}  \mathrm{d}\zeta \right) \bigg|\leq 9\varepsilon_0 e^{-x_2}.
\end{align} 
Recall that $x_1=x-R_f$ and $x_2=x$. Combining \eqref{Step74}, \eqref{Step75} and \eqref{Step76}, we conclude that 
\begin{align}\label{Step77}
	& \bigg| \E \left(1- e^{-\sum_{|u|=n} f(V(u)- (\alpha_n -x))} 1_{\left\{ M_n> \alpha_n -x \right\}}\right)- e^{-x }  C^*(f) \bigg| \leq 18\varepsilon_0 e^{-x},
\end{align}
where 
\begin{align}\label{Step78}
	C^*(f)&= e^{R_f}\ell_\infty m^a    \sum_{j=0}^\infty \E  \left(  e^{-M_j}  \int_{0} ^{\infty} e^{-\zeta}\left(1- e^{-\sum_{|v|=j} f(V(v)-M_j-\zeta+R_f)}  \right)\mathrm{d}\zeta\right) \nonumber\\
	&\quad +\ell_\infty m^a    \sum_{j=0}^\infty \E  \left(  e^{-M_j}  \int_{0} ^{\infty} e^{-\zeta}e^{-\sum_{|v|=j} f(V(v)-M_j-\zeta)}  \mathrm{d}\zeta\right).
\end{align}
It remains to check that the above definition of $C^*(f)$ is equal to that in \eqref{Def-of-C-star-f}. 
Define $S(\zeta):=1- e^{-\sum_{|v|=j} f(V(v)-M_j-\zeta)}$, then it is easy to see that $S(\zeta)=0$ when $\zeta\leq -R_f$ since $V(v)-M_j -\zeta\geq -\zeta$. Therefore, 
\begin{align}
	& e^{R_f}\int_0^\infty e^{-\zeta}S(\zeta-R_f)\mathrm{d}\zeta + \int_0^\infty e^{-\zeta} (1-S(\zeta))\mathrm{d}\zeta\nonumber\\
	&=\int_{\R} e^{-(\zeta-R_f)}S(\zeta-R_f)\mathrm{d}\zeta + 1- \int_0^\infty e^{-\zeta} S(\zeta)\mathrm{d}\zeta = 1+ \int_{0}^\infty  e^{\zeta} S(-\zeta)\mathrm{d}\zeta.
\end{align}
Putting this back to \eqref{Step78} yields that 
\begin{align}
	C^*(f)& = \ell_\infty m^a    \sum_{j=0}^\infty \E  \left(  e^{-M_j}   \left(1+ \int_0^\infty e^{\zeta} \left(1- e^{-\sum_{|v|=j} f(V(v)-M_j+\zeta)} \right) \mathrm{d} \zeta\right)\right), 
\end{align}
which implies the desired result.


\hfill$\Box$

\bigskip
\noindent

\bigskip
\noindent

\newpage

\begin{singlespace}
	
\end{singlespace}

\vskip 0.2truein
\vskip 0.2truein

\noindent{\bf Xinxin Chen:} School of Mathematical Sciences, Beijing Normal University,  Beijing, 100871, P.R. China. Email: {\texttt xinxin.chen@bnu.edu.cn}

\smallskip

\noindent{\bf Haojie Hou:}  School of Mathematics and Statistics, Beijing Institute of Technology,   Beijing 100081, P.R. China. Email: {\texttt houhaojie@bit.edu.cn}


\small
	\begin{thebibliography}{99}
		
		\bibitem{AR2009} Addario-Berry, L. and Reed, B.: Minima in branching random walks. \emph{Ann. Probab.} \textbf{37} (2009) 1044--1079.
		
		\bibitem{Aidekon2013} A\"id\'ekon, E.: Convergence in law of the minimum of a branching random walk. \emph{Ann. Probab.} \textbf{41} (2013) 1362--1426. 
		
		\bibitem{ABBS} A\"id\'ekon, E., Berestycki ,J., Brunet, \'{E}. and Shi, Z.: Branching Brownian motion seen from its tip. \emph{Probab. Theory Relat. Fields.} \textbf{157}(2013) 405--451.
		
		\bibitem{AHS} A\"id\'ekon, E., Hu, Y. and Shi, Z.: Boundedness of discounted tree sums. arXiv: 2409.01048.
		
		\bibitem{ABK} Arguin, L.-P., Bovier, A. and Kistler, N.: The extremal process of branching Brownian motion. \emph{Probab. Theory Relat. Fields.} \textbf{157}(2013) 535--574.
		
		\bibitem{BHM2018} Barral, J., Hu, Y. and Madaule, T.: The minimum of a branching random walk outside the boundary case. \emph{Bernoulli} \textbf{24}(2) (2018) 801--841.
		
		\bibitem{Biggins1976} Biggins, J.D.: The first- and last-birth problems for a multitype age-dependent branching process. \emph{Adv. in Appl. Probab.} \textbf{8} (1976) 446--459.
		
		\bibitem{BK2005}
		Biggins, J.D. and Kyprianou, A.E.: Fixed points of the smoothing transform: The boundary case. \emph{Electron. J. Probab.} \textbf{10} (2005) 609--631.
		
		\bibitem{Bramson} Bramson, M. D.: 
		Minimal displacement of branching random walk.
		\emph{Z. Wahrsch. Verw. Gebiete} \textbf{45} (2) (1978) 89--108.
		
		\bibitem{DL1983} Durrett, R. and Liggett, T.M.: Fixed points of the smoothing transformation. \emph{Z. Wahrsch. Verw. Gebiete} \textbf{64} (1983) 275--301.
		
		\bibitem{Gut2009} Gut, A.: \emph{Stopped Random walks: Limit Theorems and Applications}, 2nd ed. \emph{Springer Series in Operations Research and Financial Engineering}. New York: Springer. 
		
		
		\bibitem{Hammersley1974} Hammersley, J. M.: Postulates for subadditive processes. \emph{Ann. Probab.} \textbf{2} (1974) 652--680. 
		
		\bibitem{HuShi2009} Hu, Y. and Shi, Z.: Minimal position and critical martingale convergence in branching random walks, and directed polymers on disordered trees. \emph{Ann. Probab.} \textbf{37} (2009) 742--789.
		
		\bibitem{Jaffuel} Jaffuel, B.: The critical barrier for the survival of the branching random walk with absorption. \emph{Ann. Inst. H. Poincar\'{e}  Probab. Stat.} \textbf{48} (2012) 989--1009.
		
		\bibitem{KaPe1976} Kahane, J.-P. and Peyri\`ere, J.: Sur certaines martingales de Benoit Mandelbrot. \emph{Adv. in Math.} \textbf{22} (1976) 131--145.
		
		\bibitem{Kingman1975} Kingman, J.F.C.: The first birth problem for an age-dependent branching process. \emph{Ann. Probab.} \textbf{3} (1975) 790--801.
		
		\bibitem{Ky2000} Kyprianou, A.E.: Martingale convergence and the stopped branching random walk. \emph{Probab. Theory Related Fields} \textbf{116} (2000) 405--419. 
		
		\bibitem{Lyons}
		Lyons, R.: A simple path to Biggins' martingale convergence for branching random walk. In
		\emph{Classical and Modern Branching Processes}. IMA Vol. Math. Appl.
		\textbf{84} (1997), 217--221. Springer, New York.
		
		\bibitem{Madaule} Madaule, T.: Convergence in law for the branching random walk seen from its tip. \emph{J. Theoret. Probab.} \textbf{30} (2017) 27--63.
		
		\bibitem{Petrov1975} Petrov, V. V.:
		Sums of independent random variables.
		Springer-Verlag, New York-Heidelberg, 1975. x+346 pp.
		
		
	\end{thebibliography}
\end{document}